\documentclass[10pt, reqno]{article}
\usepackage{hyperref}
\usepackage{a4wide}
\usepackage{amsmath}
\usepackage{amsthm}
\usepackage{amssymb}
\usepackage{amsfonts}
\usepackage{color}
\usepackage{enumerate}
\usepackage{graphicx}
\usepackage{graphics}
\usepackage[french,english]{babel}
\usepackage{amsmath}
\usepackage{a4wide}
\usepackage{amssymb}
\usepackage{amsfonts}

 \usepackage{epsfig}
 \usepackage{pspicture}
  \usepackage{pstricks}
  \usepackage{eepic}

\usepackage{amsfonts}
\usepackage{color}
\usepackage{enumerate}
\def\n{\nabla}
\def\pr{\textrm{pr}(x)}
\def\l{\lambda}
\def\la{\lambda}
\def\p{\partial}
\def\D{\Delta}
\def\bou{\frac{\phi}{r_0^\l}}
\def\caat{\frac{1}{r_{0}^{\lambda}}}
\def\catt{\left(\frac{|x|}{r_0}\right)^{\lambda}}
\def\cat2{\left(\frac{1}{r_0}\right)^{2\lambda}}
\def\cat3{\left(\frac{|x|}{r_0}\right)^{2\lambda}}
\def\rr{\mathbb{R}}
\def\ccat{\left(\frac{|x|}{r_0}\right)^{\lambda-2}}

\def\dd{d\theta_{N-1}\ldots d\theta_2d\theta_1 dr}

\def\into{\int_{\Omega}}
\def\hoi{H_{0}^{1}(\Omega)}

\def\dd{\mathrm{d}}

\def\q{\quad}

\def\eps{\varepsilon}

\def\eps{\varepsilon}

\def\into{\int_{\Omega}}
\def\hoi{H_{0}^{1}(\Omega)}

\newcommand{\intos}{\int_{\tilde{\Omega}_{r_1}}}
\def\eps{\varepsilon}

\def\be{\begin{equation}}
\def\ee{\end{equation}}

\def\Omr{\mathcal{O}}
\def\Oms{\tilde{\mathcal{O}}}
\def\Omg{\tilde{\Omega}_{r_1}}
\def\iiQ{\iint\limits_{\Omega\times(0,T)}}
\def\iir0{\iint\limits_{\tilde{\Omega}_{r_0}\times(0,T)}}

\def\iiq{\iint\limits_{\omega_0\times(0,T)}}

\def\iiQr{\iint\limits_{\Omr\times(0,T)}}
\def\iiQs{\iint\limits_{\Oms\times(0,T)}}
\def\is{\iint\limits_{\Gamma \times(0,T)}}

\def\dx{\ \mathrm{dx}}

\def\dt{\ \mathrm{dt}}
\def\ds{\ \mathrm{ds}}
\def\|{\Big |}
\def\({\Big (}
\def\){\Big )}
\def\[{\Big[}
\def\]{\Big]}
 \newtheorem{thm}{Theorem}[section]
 
 \newtheorem{lem}[thm]{Lemma}
 \newtheorem{prop}[thm]{Proposition}
 \theoremstyle{definition}
 
 \theoremstyle{remark}
 \newtheorem{rem}[thm]{Remark}
 \numberwithin{equation}{section}

\begin{document}
\begin{center}
{\Large Controllability of the  heat equation with an inverse-square
potential localized on the boundary}\\[5mm]

\bf Cristian Cazacu\footnote{Research group of the projects PN-II-ID-PCE-2011-3-0075 and PN-II-ID-PCE-2012-4-0021, ``Simion Stoilow" Institute of Mathematics of the Romanian Academy, P.O. Box 1-764, 014700 Bucharest, Romania.\\
E-mail: cristi\_cazacu2002@yahoo.com.}
\end{center}

\begin{abstract}
This article is devoted  to analyzing control properties  for the heat
equation with singular potential $-\mu/|x|^2$  arising at the
boundary of a smooth domain $\Omega\subset \rr^N$, $N\geq 1$. This
problem was firstly studied by Vancostenoble and Zuazua
\cite{heatjudith} and then generalized by Ervedoza \cite{sylvain}
 in the context of
interior singularity. Roughly speaking, these results  showed  that
for  any value of parameters $\mu\leq \mu(N):=(N-2)^2/4$,   the
corresponding parabolic system can be controlled to zero with the
control distributed in any open subset of the domain.  The critical
value $\mu(N)$ stands for the best constant in the Hardy inequality
with interior singularity.

When considering the case of boundary singularity a better critical
Hardy constant is obtained, namely $\mu_{N}:=N^2/4$.

In this article we extend the previous results  in
\cite{heatjudith}, \cite{sylvain},
 to the case of
boundary singularity. More precisely,  we show that for any $\mu
\leq \mu_N$, we can lead the system to zero state using a
distributed control in any open subset.

 We emphasize that our results cannot be obtained
straightforwardly from the previous works \cite{heatjudith},
\cite{sylvain}. 

\end{abstract}


\section{Introduction}

In this article we present some new results concerning the exact
controllability  of the heat equation with singular quadratic
potential $-\mu/|x|^2$.

From a mathematical viewpoint, the study of problems with
inverse-square potentials is motivated by models which appear for
instance in the context of combustion theory \cite{MR0110868}, \cite{MR0092663} and quantum
mechanics \cite{simon}.

Evolution problems with the potential $-\mu/|x|^2$ have been
intensively studied in the recent decades. Among such results, we remind the
pioneering work by Baras and Goldstein \cite{baras} in which they
considered the corresponding heat-like equation with the singularity
localized in the interior of a bounded domain $\Omega\subset\rr^N$,
$N\geq 2$ (If $N=1$ they deleted the origin so that $0\in \Gamma$, where $\Gamma$  denotes the boundary of $\Omega$). They
derived necessary and sufficient conditions for such systems to
be well-posed. More precisely, they showed the well-posedness
holds true whenever $\mu\leq (N-2)^2/4$, whereas if $\mu>(N-2)^2/4$
there is instantaneous blow-up for the solution. The critical
value $(N-2)^2/4$ is the best constant in the corresponding Hardy
inequality (see e.g. \cite{hard}, \cite{hardy-polya}).
 Later on, the issue of singular or degenerated potentials has also been analyzed
 by the control community. Among the pioneering related works we mainly refer to the paper by  Cannarsa, Martinez and Vancostenoble
 \cite{martinez} and references therein studying the control
 of  parabolic equations degenerating at origin.

  The authors in \cite{judith} analyzed the control
 and stabilization properties
 of the wave equation with the singular potential $-\mu/|x|^2$ located in the interior of the domain. Then they showed in
Vancostenoble and Zuazua \cite{heatjudith} that the corresponding
 heat  equation can be controlled by a distributed control which
 sourrounds the singularity. This result was generalized in
 Ervedoza \cite{sylvain} where any geometrical constraint
 of the control region was removed.
  Recently,
 the work \cite{sylvain}
 has been slightly improved in \cite{vanc} when studying some
 applications to inverse problems. In all situations above the authors
 showed that well-posedness, control and stabilization are very
 much related to the classical Hardy inequality in which the best
 constant is $(N-2)^2/4$.

In this paper we consider the heat equation with the potential
$-\mu/|x|^2$, where the singularity $x=0$ is located on the boundary
$\Gamma$ of a bounded domain  $\Omega\subset\rr^N$, $N\geq 1$. For the sake of clarity, let the boundary of $\Omega$  be smooth enough (some intermediate results require $C^4$ smoothness for the boundary $\Gamma$).

This work aims to extend to the case of boundary singularity the results of
paper \cite{sylvain} which provides the most general control results
in the case when the singularity is localized in the interior of the
domain. We point out that our results cannot be deduced straightforwardly
from the case of interior singularity and requires an independent
analysis. Our main tool relies on Carleman estimates which is
the classical way to prove observability properties for parabolic
systems. The major difficulty consists in finding proper weight
functions to develop efficient Carleman estimates.
 In our case, the weights in \cite{sylvain}, \cite{heatjudith}, are
 not even
allowed  to recover the results in the range of parameters $\mu\leq
(N-2)^2/4$ shown in the case of interior singularity. A proper
modification of the weights in \cite{sylvain}, will be done here.

Before entering into details, let us fix some ideas.

 Let $\Omega \subset \rr^N$, $N\geq 1$,  such that $0 \in \Gamma$, and let
$\omega\subset\Omega$ be a non-empty open set. Assume also that
$T>0$ is fixed. We are interested in the question of controllability
 of
the following problem \be\label{controleq}
    \left\{
        \begin{array}{lll}
            \displaystyle \partial_t u - \Delta u - \frac{\mu}{|x|^2} u = f,
            &\qquad  (x,t) \in \Omega\times (0,T), \\
            u(x,t)=0, &\qquad (x,t) \in \Gamma \times (0,T),
            \vspace{0.3ex}\\
            u(x,0 ) = u_0(x), &\qquad x \in\Omega,
        \end{array}
    \right.
\ee where $u_0\in L^2 (\Omega)$ and $f\in L^2(\Omega\times (0, T))$
is a function supported in the control region $\omega$.

The null-controllability problem reads as follows: Given any $u_0
\in L^2(\Omega)$, find a function $f\in L^2(\omega\times(0,T))$ such
that the solution of \eqref{controleq} satisfies
\be\label{nullcontrol}
    u(x,T) =0, \qquad  x \in \Omega.
\ee

In order to discuss the well-posedness and null-controllability of
\eqref{controleq} we need to establish the proper functional
framework which corresponds to the problem. The crucial role of this
issue is played by a new critical value of $\mu$ which determines
the features of system \eqref{controleq}. More precisely,
  when moving the singularity from the interior to the boundary,
   the critical Hardy constant jumps from $(N-2)^2/4$ to
the critical value
\begin{equation}\label{newcritical}
\mu_N:=\frac{N^2}{4}.
\end{equation}
 This is guaranteed by the improved   Hardy inequalities
  with boundary singularities
 stated in Propositions \ref{prop1}-\ref{prop2} as follows.
  \begin{prop}\label{prop1}
    Let $\Omega\subset \rr^N$, $N\geq 1$, be a domain satisfying  $0\in \Gamma$. Then, for any $\mu \leq  \mu_N$ and  any $0\leq \gamma<2$, 
  there exists a constant $C_1$ depending on $\gamma, \mu$ and  $\Omega$,  such that  the inequality
\be\label{Hardy}
    \qquad     \int_{\Omega} \frac{u^2}{|x|^\gamma}\dx
    + \mu \into
    \frac{u^2}{|x|^2} \dx \leq \into |\nabla u|^2 \dx + C_1\into u^2\dx, \ee holds for all $u\in \hoi$.
\end{prop}
\begin{prop}\label{prop2}
    Let $\Omega\subset \rr^N$, $N\geq 1$, be a domain satisfying  $0\in \Gamma$. Then, for any $\mu
\leq \mu_N$ and any $0\leq  \gamma<2$, there exist two constants $C_2$ and  $C_3$ depending on $\gamma, \mu$ and  $ \Omega$,   such that
\begin{equation}\label{imprhardy}
C_2 \into u^2 \dx + \into \left( |\n u|^2 -\mu\frac{u^2}{|x|^2} \right)
\dx\geq C_3 \into \left(|x|^{2-\gamma} |\n u|^2
+\frac{u^2}{|x|^\gamma}\right)\dx, \q \forall u\in \hoi.
\end{equation}
\end{prop}
These results will be used in the proof of the Carleman estimates.

Here we skip the proof of  Proposition \ref{prop1}  since it is a direct consequence of the  inequalities with boundary singularities and logarithmic reminder terms showed in   \cite{cristiCRAS}, \cite{musina},
\cite{Fall}.

On the other hand,  Proposition \ref{prop2} was proved in   \cite{cristisumJFA} (see Theorem 1.1) in the particular case  $\gamma=0$.  Following the proof in \cite{cristisumJFA}, in the Appendix we give a rigorous justification of inequality \eqref{imprhardy} for any $\gamma\in [0, 2)$.

Next we can formulate the main result of this paper.
\begin{thm}[Null-Controllability]\label{boule}
    Let $\Omega\subset \rr^N$, $N\geq 1$, be a domain satisfying  $0\in \Gamma$ and assume $\mu\leq \mu_N$. Given any non-empty
open set $\omega\subset \Omega$, for any time $T>0$ and any initial
data $u_0\in L^2(\Omega)$, there exists a control $f\in
L^2(\omega\times (0, T))$ such that the solution of
\eqref{controleq} satisfies \eqref{nullcontrol}.
\end{thm}

\begin{rem}
The authors in \cite{heatjudith} proved the null-controllability of system \eqref{controleq} with interior singularity,  acting with a control supported in an annulus surrounding the origin. They derived their result by means of the spherical harmonics decomposition,  reducing the problem to the one-dimensional case in which the singularity arises at the origin and the control $\omega$ is distributed in an interval, say, $\Omega=(0,1)$.
In other words, the result obtained in \cite{heatjudith} is equivalent to the result of Theorem \ref{boule} in the case $N=1$. Of course, in this paper we are concerned about the validity of Theorem \ref{boule} in the non-trivial case $N\geq 2$.
\end{rem}
\vspace{0.2cm}

\begin{figure}[h]
\begin{center}
\scalebox{1} 
{
\begin{pspicture}(0,-1.3430756)(9.205419,1.3430754)
\psbezier[linewidth=0.04](0.5014236,-0.45739147)(0.0,0.40781045)(0.68648326,0.39901936)(1.6014236,0.8026085)(2.5163639,1.2061977)(2.573247,1.3230755)(3.4014237,0.7626085)(4.2296004,0.20214146)(4.6694517,0.12829262)(3.9814236,-0.5973915)(3.2933955,-1.3230755)(3.4985194,-0.1335477)(2.5014236,-0.05739148)(1.5043277,0.018764732)(1.0028472,-1.3225935)(0.5014236,-0.45739147)
\psellipse[linewidth=0.04,dimen=outer](1.7714236,0.5826085)(0.25,0.12)
\psline[linewidth=0.04cm](1.5814236,0.5426085)(1.7814236,0.6426085)
\psline[linewidth=0.04cm](1.8014235,0.6626085)(1.8414236,0.68260854)
\psline[linewidth=0.04cm](1.7414236,0.50260854)(1.9214236,0.62260854)
\psbezier[linewidth=0.04](6.2214236,-0.15739149)(6.646007,-1.0627803)(8.777428,-0.83636326)(8.981423,0.14260852)(9.185419,1.1215802)(8.433647,1.1470783)(7.4414234,1.0226085)(6.4492,0.8981387)(5.79684,0.7479973)(6.2214236,-0.15739149)
\psellipse[linewidth=0.04,dimen=outer](8.231423,0.48260853)(0.25,0.12)
\psline[linewidth=0.04cm](8.041424,0.4426085)(8.241424,0.5426085)
\psline[linewidth=0.04cm](8.261423,0.56260854)(8.301424,0.5826085)
\psline[linewidth=0.04cm](8.201424,0.4026085)(8.381424,0.5226085)
\psdots[dotsize=0.172](2.3814235,-0.05739148)
\psdots[dotsize=0.172](7.5814238,-0.7173915)
\usefont{T1}{ptm}{m}{n}
\rput(2.9,0.6){$\Omega$}
\rput(2.5,-0.4){$x=0$}
\rput(1.3,0.5){$\omega$}
\rput(7, 0.5){$\Omega$}
\rput(7.8,-0.95){$x=0$}
\rput(8.4,0.15){$\omega$}
\end{pspicture}
}\caption{The control region $\omega$ may be any open subset  independent on the local shape of $\Omega$ at the singularity $x=0$.}\label{region}
\end{center}
\end{figure}
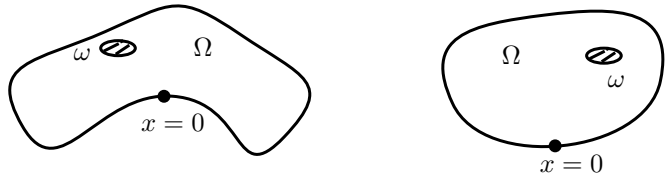

Theorem \ref{boule}  says that our main  control results  do not require  any constraints for the control region $\omega$ in the sense that $\omega$ is allowed to be any open subset of $\Omega$, no matter what the geometry of $\Omega$ is,  as
 depicted in Figure \ref{region}.
\vspace{0.4cm}

 Following
the by now classical HUM method (cf. \cite{lions1}), the
controllability property is equivalent to an observability
inequality for the adjoint system \be\label{adjoint}
    \left\{
        \begin{array}{lll}
            \displaystyle \partial_t w + \Delta w + \frac{\mu}{|x|^2} w = 0, &\qquad & (x,t)
            \in \Omega\times (0,T), \\
            w(x,t) =0, &\qquad & (x,t) \in \Gamma \times (0,T), \vspace{0.3ex}\\
            w(x,T ) = w_T(x), &\qquad & x \in\Omega.
        \end{array}
    \right.
\ee More precisely, when $\mu\leq \mu_N$, we need to prove that
there exists a constant $C_T$ such that for all $w_T \in L^2(\Omega)$,
the solution of \eqref{adjoint} satisfies \be\label{obs}
    \int\limits_{\Omega} |w(x,0)|^2\dx \leq C_T \iint\limits_{\omega\times(0,T)}
     |w(x,t)|^2 \dx \dt.
\ee

In order to prove \eqref{obs}, we will use a particular Carleman
estimate, which is by now a classical technique in control theory. 
 Indeed, the Carleman estimate we will derive later implies that for
any solution $w$ of \eqref{adjoint},
 \be\label{obsa}
    \iint\limits_{\Omega\times(\frac{T}{4},\frac{3T}{4})} |w(x,t)|^2
     \dx \dt \leq C_T \iint\limits_{\omega\times(0,T)} |w(x,t)|^2
     \dx\dt.
\ee
Let us show that \eqref{obsa} implies \eqref{obs}. Indeed,
multiplying the system \eqref{adjoint} by $w$ and integrating over
$\Omega$ we formally obtain

$$ \frac{1}{2}\frac{\dd}{\dt}\int_{\Omega}
w^2(x,t)\dx=\int_{\Omega} |\nabla w(x, t)|^2\dx-\mu
\int_{\Omega}\frac{w^2(x, t)}{|x|^2}\dx.$$

From \eqref{Hardy} we have that
$$\frac{1}{2}\frac{\dd}{\dt}\int_{\Omega}
w^2(x,t)\dx\geq -c\int_{\Omega}w^2(x,t)\dx,
$$ for some constant $c>0$ depending on $C_1$ in \eqref{Hardy}.
  Then  we get that the
function $t\mapsto e^{2ct}||w(\cdot,
t)||_{L^2(\Omega)}^{2}$ is increasing and  we have

\be\label{eq2} \int_{T/4}^{3T/4}\into w^2(x,t)\dx \dt \geq
\frac{T}{2}e^{-3Tc/2}\into w^2(x,0)\dx.
 \ee
From here and  \eqref{obsa} we obtain \eqref{obs}.

\hspace{0.8cm}

%

\noindent {\bf Well-posedness via Hardy inequality}
%

Let us fix $0\leq \gamma<2$ and define the set
\begin{equation}\label{set}
 \mathcal{L}^\gamma:= \Bigg \{ C\geq 0 \ \textrm{ s. t. } \ \inf_{u\in \hoi}
 \frac{\into \left(|\n u|^2-\mu_Nu^2/|x|^2+C u^2 \right) \dx
 }{\into
 u^2/|x|^\gamma \dx} \geq 1 \ \Bigg\}.
\end{equation}

Of course, $\mathcal{L}^\gamma$ is non empty since from  inequality
\eqref{Hardy} we have that $|C_1|\in \mathcal{L}^\gamma$. Next we define
\begin{equation}\label{optimalconstant}
\mathcal{C}_{0}^{\gamma}=\inf_{C\in \mathcal{L}^\gamma}C.
\end{equation}

Then, for any  $\mu\leq \mu_N$ we introduce the Hardy
functional
\begin{equation}\label{eq130}
B_{\mu}^{\gamma}(u):=\int_{\Omega} |\nabla u|^2 \dx-\mu\int_{\Omega}
\frac{u^2}{|x|^2}\dx+\mathcal{C}_{0}^{\gamma}\int_{\Omega} u^2 \dx,
\end{equation}
which is positive for any test function due to inequality
\eqref{Hardy} and the election of $\mathcal{C}_0^\gamma$. Then we define
the  corresponding Hilbert space $H_{\mu}^{\gamma}$ as the closure of
$C_{0}^{\infty}(\Omega)$ in the norm induced by  $B_{\mu}^{\gamma}(\cdot)$.

By the definition of $\mathcal{C}_{0}^{\gamma}$, if $\mu\leq  \mu_N$ we obtain
\begin{multline}\label{normsubcritical}
\left(1- \frac{\mu^+}{\mu_N}\right)\into \left(|\n u|^2 +\mathcal{C}_{0}^{\gamma}u^2\right) \dx +\frac{\mu^+}{\mu_N}\into \frac{u^2}{|x|^\gamma} \dx \\
\leq||u||_{H_{\mu}^{\gamma}}^2 \leq \\ \left(1+\frac{\mu^-}{\mu_N}\right)\into \left(|\n u|^2+\mathcal{C}_{0}^{\gamma} u^2\right) \dx,
\end{multline}
where $\mu^+:=\max\{0, \mu \}$ and $\mu^{-}:=\max \{0, -\mu\}$.

Observe that for any $\mu< \mu_N$ the identification $H_{\mu}^{\gamma} =\hoi$
holds true due to the equivalence of the corresponding norms in \eqref{normsubcritical}.
Besides, in the critical case $\mu=\mu_N$, $H_{\mu_N}^{\gamma}$ is slightly
larger than $\hoi$. However,  using cut-off arguments near the
singularity (see e.g. \cite{VazZua}) we can show that
\begin{equation}\label{slightly}
||u||_{H_{\mu_N}^{\gamma}}\geq C_{\eps}^{\gamma} ||u||_{H^{1}(\Omega \setminus
\overline{B(0, \eps)})}, \q \forall  u\in H_{\mu_N}^{\gamma},
\end{equation}
where $C_{\eps}^{\gamma}$ is a constant going to zero as $\eps$ tends to zero
and $\overline{B(0, \eps)}$ denotes the closure of the ball of radius
$\eps$ centered at the origin.

Let us define now the operator $A_\mu^\gamma:=-\D-\mu/|x|^2+\mathcal{C}_0^\gamma
I$ together with its domain as
\begin{equation}\label{domain}
D(A_{\mu}^{\gamma}):=\{ u \in H_{\mu}^{\gamma} \ | \  A_{\mu}^{\gamma} u \in L^2(\Omega)\}.
\end{equation}
The norm of this operator is given by
\begin{equation}\label{operator}
||u||_{D(A_{\mu}^{\gamma})}=||u||_{L^2(\Omega)}+ ||A_{\mu}^{\gamma} u||_{L^2(\Omega)}\q
\forall \mu\leq \mu_N.
\end{equation}
With these definitions, by standard semigroup-theory one can show
that for any $\mu\leq \mu_N$ the operator $(A_{\mu}^{\gamma}, D(A_{\mu}^{\gamma}))$
generates an analytic semigroup in the pivot space $L^2(\Omega)$ for
the equation \eqref{controleq}. For more details we refer to Hille-Yosida theory in
Chapter 7, p. 190, \cite{brezis}, which can be adapted in the context of the
space $H_{\mu}^{\gamma}$ introduced above. For our particular problem the interested reader can also consult  Theorem  II.1, p. 3, \cite{vanc} which refers to well-posedness in the case of an interior singularity.

\bigskip

\vspace{0.1cm}

First let us briefly present the structure of the paper. In Section \ref{preuvesanscalcul} we design new weights of the Carleman estimates which are adapted to our problem. The corresponding Carleman inequality stated in Theorem \ref{teoremnewvariable} leads to the controllability result of Theorem \ref{boule} as described above.   The proof of Theorem \ref{teoremnewvariable} is  obtained by gluing the inequalities of Lemmas \ref{pozbd}-\ref{Ir}. Due to technical difficulties, the main body of the paper is devoted to proving these lemmas. In view of that,  in Section \ref{babizduca} we develop some preliminary technical results to conclude with the proofs of Lemmas \ref{pozbd}-\ref{Ir} in Section \ref{desteptutz}.  Finally, we end up with an Appendix (cf. Section \ref{7sec}) where we show the validity of the Hardy-type inequality announced in Proposition \ref{prop2} above.

\section{Null controllability in the case $\mu\leq N^2/4$} \label{preuvesanscalcul}

First of all, to simplify the presentation, we assume that $0 \notin
\bar \omega$ otherwise it is straightforward since the control acts
locally near the singularity. We also assume that $\Omega\cap
\overline{ B_1(0)} \cap \bar \omega$ is empty. This can always be done by a scaling
argument. In the sequel we also consider a nonempty subset
$\omega_0\subset\subset \omega$ whose role will be emphasized in the next sections.

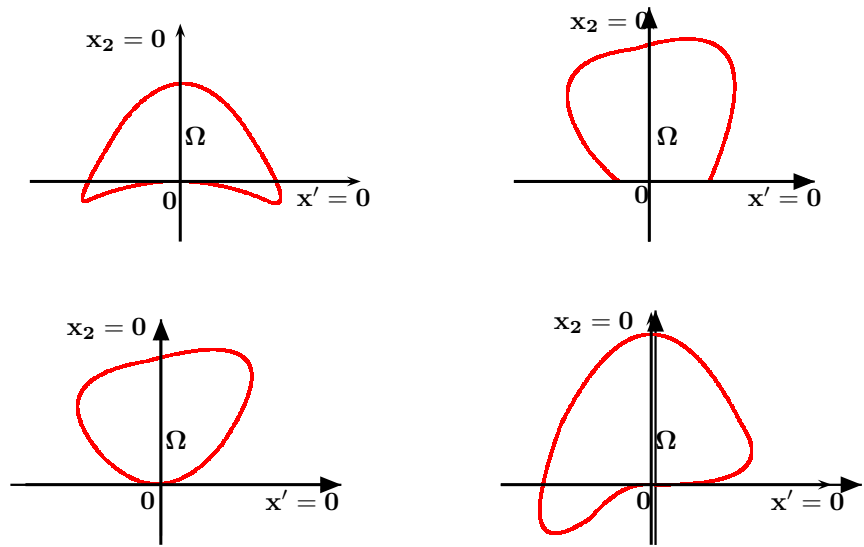
\begin{figure}[h!]
\begin{center}
\setlength{\unitlength}{0.4cm}
\begin{picture}(15,10)
\linethickness{0.3mm}
   \textcolor{red}{\qbezier(4,2.4)(7,3.65)(10,2.3)
   \qbezier(4,2.4)(3.15,1.9)(5.1,5)
   \qbezier(10,2.3)(11,2)(9,5)
   \qbezier(5.1,5)(7.1,7.5)(9,5)}
   \put(6.4,2.1){\bf{0}}
   \put(7.5,4.5){\makebox(0,0){\emph{$\mathbf{\Omega}$}}}
   \put(12.1,2.5){\makebox(0,0){$\mathbf{x'=0}$}}
   \put(5.2,7.7){\makebox(0,0){$\mathbf{x_2=\bf{0}}$}}
   \psline[linewidth=0.04cm,arrowsize=0.05291667cm
2.0,arrowlength=1.4, arrowinset=0.4]{->}(0.8,1.2)(5.2,1.2)
\psline[linewidth=0.04cm,arrowsize=0.05291667cm
2.0,arrowlength=1.4,arrowinset=0.4]{->}(2.8, 0.4)(2.8, 3.3)
  \end{picture}
\setlength{\unitlength}{0.4cm}
 \begin{picture}(14,10)
\linethickness{0.3mm}
 \textcolor[rgb]{0.00,0.00,0.50}{ \put(2.5,3){\vector(1,0){10}}}
  \put(7,1){\vector(0,1){7.8}}
   \textcolor[rgb]{1.00,0.00,0.00}{\qbezier(5,4)(5.8,3.05)(6,3)
   \qbezier(5,4)(3,7)(6.5,7.4)
   \qbezier(6.5,7.4)(11.5,9)(9,3)
 } 
   \put(7.6,4.5){\makebox(0,0){\emph{$\mathbf{\Omega}$}}}
   \put(11.5,2.5){\makebox(0,0){$\mathbf{x'=0}$}}
   \put(5.7,8.3){\makebox(0,0){$\mathbf{x_2=\bf{0}}$}}
   \put(6.5,2.3){\bf{0}}
   \psline[linewidth=0.04cm,arrowsize=0.05291667cm
2.0,arrowlength=1.4, arrowinset=0.4]{->}(1,1.2)(5,1.2)
\psline[linewidth=0.04cm,arrowsize=0.05291667cm
2.0,arrowlength=1.4,arrowinset=0.4]{->}(2.8, 0.4)(2.8, 3.4)
  \end{picture}\\
  \begin{picture}(16,7)
\linethickness{0.3mm}
  \put(2,3){\vector(1,0){11}}
  \put(7,1.3){\vector(0,1){7.2}}
   \textcolor{red}{\qbezier(5,4)(7,1.9)(9,4.3)
   \qbezier(5,4)(2.9,6.5)(6.5,7.1)
   \qbezier(6.5,7.1)(11.9,8.6)(9,4.3)}
   \put(7.5,4.5){\makebox(0,0){\emph{$\mathbf{\Omega}$}}}
   \put(11.7,2.5){\makebox(0,0){$\mathbf{x'=0}$}}
   \put(5.2,8.2){\makebox(0,0){$\mathbf{x_2=\bf{0}}$}}
    \put(6.3,2.2){\bf{0}}
   \psline[linewidth=0.04cm,arrowsize=0.05291667cm
2.0,arrowlength=1.4, arrowinset=0.4]{->}(1,1.2)(5.1,1.2)
\psline[linewidth=0.04cm,arrowsize=0.05291667cm
2.0,arrowlength=1.4,arrowinset=0.4]{->}(2.8, 0.4)(2.8, 3.3)
    \end{picture}
  \begin{picture}(14,10)
\linethickness{0.3mm}
  \put(2,3){\vector(1,0){12}}
  \put(7.15,1){\vector(0,1){7.8}}
  \textcolor{red}{
 \qbezier(7,3)(11.4, 2.9)(10, 5)
   \qbezier(5,1.9)(6,3.1)(7,3)
   \qbezier(5,1.9)(2.2,0)(4,5)
    \qbezier(4,5)(7,11)(10,5)
   }
   \put(7.5,4.5){\makebox(0,0){\emph{$\mathbf{\Omega}$}}}
   \put(6.5,2.2){\bf{0}}
   \put(12.2,2.5){\makebox(0,0){$\mathbf{x'=0}$}}
   \put(5.1,8.4){\makebox(0,0){$\mathbf{x_2=\bf{0}}$}}
   \psline[linewidth=0.04cm,arrowsize=0.05291667cm
2.0,arrowlength=1.4, arrowinset=0.4]{->}(0.8,1.2)(5.2,1.2)
\psline[linewidth=0.04cm,arrowsize=0.05291667cm
2.0,arrowlength=1.4,arrowinset=0.4]{->}(2.8, 0.4)(2.8, 3.5)
  \end{picture}\quad
  \caption{\label{ffig1} The geometry of a domain $\Omega\subset \rr^2$  containing the origin on the boundary and verifying $\vec{n}(0)=-e_2$.  Such a two-dimensional domain  might have four configurations (convex, concave, flat or changing convexity at the origin). The condition \eqref{impp} is reflected in the fact that the tangent at $x=0$ to $\Gamma$ coincides to the $x'$-axis. In particular, such a domain satisfies the condition \eqref{impcond1}.}
   \end{center}
\end{figure}

In addition, without loss of generality, since the operator $\p_t-\Delta-\mu/|x|^2$ is invariant under rotations centered at $x=0$, we may assume that
\begin{equation}\label{impp}
\vec{n}(0) =-e_N,
\end{equation}
where $\vec{n}(0)$ is the outward normal vector to $\Gamma$ at $x=0$ and $e_N$ is the $N$-th unit vector of the canonical basis $\{e_i\}_{i=1, \ldots, N}$ in $\rr^N$ (in Figure \ref{ffig1} we emphasize  the condition \eqref{impp} for a domain $\Omega\subset \rr^2$).
Moreover, for $N\geq 2$,  we consider  the notation
$$ x=(x', x_N) \in \rr^N, \textrm{with }  x'\in  \rr^{N-1}, x_N \in \rr. $$

On the other hand, $\Gamma$ is an $(N-1)$ - Riemannian submanifold of $\rr^N$ and admits a local smooth parametrization  $\Phi:\rr^{N-1}\rightarrow \rr$, i.e.  there exists a neighborhood $\mathcal{V}$ of $x=0$,   such that
 \begin{equation}\label{condphi}
 x_N=\Phi(x'), \quad \forall x\in \Gamma \cap \mathcal{V}.
 \end{equation}
 Combining this with \eqref{impp} we obtain that $\Phi$ has a  quadratic degeneracy as $x\rightarrow 0$, that is
 \begin{equation}\label{gradphi}
 \Phi(0)=0, \quad \n_{x'} \Phi(0)=0.
 \end{equation}
 For those reasons,  by Taylor expansion, the local parametrization of  $\Gamma$ verifies $x_N=\Phi(x')=O(|x'|^2)$, as $x\rightarrow 0$.
 As a consequence, we deduce that the points on the boundary $\Gamma$ satisfy
\begin{equation}\label{impcond1}
|x\cdot \vec{n}(x)|\leq C_\Omega |x|^2, \quad \forall x\in \Gamma,
\end{equation}
 where $C_\Omega$ is a positive constant depending on $\Omega$ and $\vec{n}(x)$ stands for the outward unit normal vector at any point  $x\in \Gamma$.

To simplify, we conclude that under the assumption \eqref{impp} (up to a rotation, this could always  be assumed), inequality \eqref{impcond1} is verified and will play a crucial role in our Carleman estimates.

%
%
%

In what follows our aim is to justify the result of Theorem
\ref{boule}. For that we will apply Carleman estimates by modifying the standard weights in \cite{fursikov} and \cite{sylvain}. 



\subsection{Carleman estimates. The choice of a proper weight}\label{CarlSub}

As already stated in Introduction, the main tool we use to address the
 observability inequality \eqref{obsa} is a Carleman estimate.

 The major problem when designing a Carleman estimate
is the choice of a smooth weight function $\sigma$, which is in
general assumed to be positive, and to blow up as $t$ goes to zero
and as $t$ goes to $T$. Hence we are looking for a weight function
$\sigma$ that satisfies:
 \begin{equation}\label{hypsigma}
    \left\{ \begin{array}{ll}
        \sigma(t,x) > 0, \quad (t,x) \in (0,T)\times \Omega,\vspace{1ex} \\
        \lim\limits_{t\rightarrow 0^+} \sigma (t,x) =\lim\limits_{t\rightarrow T^-}
        \sigma (t,x) = +\infty, \quad x \in \Omega.
    \end{array}\right.
\end{equation}
When shifting the singularity from the interior up to the boundary the
weight in Ervedoza \cite{sylvain} violates some necessary conditions
to prove the Carleman estimates; in particular, the weight $\psi$ in \cite{sylvain} blows up at
origin and this violates the fact that $\psi$ is constant on the
boundary. In the next section we design a new $\psi$ which fits our problem.

\subsubsection{ The definition of $\psi$}\label{papat}

 In what follows we need to recall some classical properties of the distance function to the boundary, say $\rho$,  which turns out to be very important and useful in order to introduce  $\psi$ and to develop our observability inequality.\\

\noindent \textsl{1. The distance function $\rho(x)$}.

As emphasized in \cite{marcus},  for $\beta>0$ let
 \begin{equation}
 \Omega_\beta=\{ x\in \Omega; \ \rho(x)<\beta \}, \quad \Sigma_\beta=\{ x\in \Omega; \ \rho(x)=\beta \},
 \end{equation}
 where $\rho(x)=\textrm{dist}(x, \Gamma)$ denotes the Euclidian distance to the
boundary.  For $\beta_0$ small enough there exists a unique point $\textrm{pr}(x)\in \Gamma$ such that
\begin{equation}\label{puffy}
\rho(x)=|x-\pr|, \quad \forall x\in \Omega_\beta, \quad \forall \beta\leq \beta_0.
\end{equation}
Moreover, the mapping
$$\Pi :\Omega_{\beta} \rightarrow (0, \beta)\times \Gamma, \quad \Pi(x)=(\rho(x), \pr)$$
is a diffeomorphism  and its inverse is given by
$$ (0, \beta)\times \Gamma \rightarrow \Omega_{\beta}, \quad  \Pi^{-1}(t, \pr)= \pr -t \ \vec{n}(\pr), \qquad \forall \beta\leq \beta_0. $$
In particular, the distance function $\rho(x)$ is smooth on $\Omega_{\beta_0}$ and satisfies both
\begin{equation}\label{turca}
\n \rho(x) = -\vec{n}(\pr), \quad \forall x \in \Omega_{\beta_0},
\end{equation}
and the Eikonal equation
 \begin{equation}\label{eikeq}
 |\n \rho(x)|=1, \quad \forall x\in \Omega_{\beta_0}.
 \end{equation}

\noindent \textsl{2. Definition of $\psi$ via $\psi_1$}

 Next we introduce a smooth function
 $\psi_1$ (at least $C^4(\overline{\Omega})$) satisfying the conditions
\begin{equation}\label{eeee}
\left\{\begin{array}{ll}
\psi_1(x)=\rho(x), & \forall x\in \Omega_{\beta_0},\\
 \psi_1(x) >\beta_0, & \forall  x \in  \Omega\setminus \overline{\Omega_{\beta_0}},\\
   \psi_1(x) \equiv \beta_0, & \forall   x\in \Sigma_{\beta_0},\\
 |\nabla \psi_1(x)|\geq \delta_0>0 & \forall  x\in \Omega \setminus
\overline{\omega_0},\end{array}\right.
\end{equation}
Such a function exists but its construction is not trivial.
Indeed, there exists a smooth positive function which extends the distance $\rho(x)$ from $\Omega_{\beta_0}$ to $\Omega\setminus \Omega_{\beta_0}$ since this feature is generally true. This function vanishes on the boundary $\Gamma$ and,  according to classical arguments of Morse theory (cf.  p. 80, \cite{coron}),  since  $\rho(x)$ satisfies the equation \eqref{eikeq} in $\Omega_{\beta_0}$, it has all (finitely many) critical points located in $\Omega\setminus \overline{\Omega_{\beta_0}}$ . Then we consider such a function and, following the construction in \cite{fursikov}  (through a diffeomorphism transformation), we move the critical points into  $\omega_0$  without modifying the function in $\Omega_{\beta_0}$. Afterwards, we obtain the existence of $\psi_1$ as in \eqref{eeee}.

To conclude, we design the function $\psi\in C^4(\overline{\Omega})$ given by
\begin{equation}\label{defimp}
\psi:=\delta (\psi_1+1),
\end{equation}
 such that $\delta \delta_0> 2 C_\Omega$ where $\delta_0$ stands for the constant in \eqref{eeee} and $C_\Omega$ is as in  \eqref{impcond1}. In particular,
 under these conditions $\psi$ satisfies  the following useful properties necessary next in the paper:
\begin{equation}\label{equ1}
\left\{\begin{array}{ll}
\psi(x)=1, & \forall x\in \Gamma,\\
 \psi(x) >1, & \forall  x \in  \Omega,\\
  |\nabla \psi(x)|> 2C_\Omega & \forall x\in \Omega \setminus
\overline{\omega_0}.\end{array}\right.
\end{equation}
Moreover, due to technical computations which will be expressed later, we fix $\delta$ such that
\begin{align}\label{bosgo}
\delta & \geq  \max \Bigg\{1, \frac{2C_\Omega}{\delta_0}, \frac{24D_{\Omega, \psi_1}R_\Omega^2}{\delta_0^2}, \frac{2}{\delta_0}, \frac{1}{\delta_0^2} \left(1+ 4 D_{\Omega, \psi_1} + |D\psi_1|_{\infty}+ 2|D^2 \psi_1|_{\infty}\right), \Bigg\},
\end{align}
where $D_{\Omega, \psi_1}$ is a constant depending only on $\Omega$ and $\psi_1$ according to Lemma \ref{crucial}, and $R_\Omega=\sup_{x\in \overline{\Omega}}|x|$. \\

\noindent {\textsl{\bf Notations:}}

\noindent Throughout the paper,  formally, for a given function $f$ we understand\\

\noindent $|f|_\infty=|f|_{L^\infty(\Omega)}$, \\
$|D f|_{\infty}= |\n f|_{L^\infty(\Omega)} $,\\
$D^2 f (\xi, \xi)= \sum_{i, j=1}^{N} \partial_{x_ix_j}^{2} f \xi_i\xi_j$, \quad $\forall \xi \in \rr^N$,   \\
$|D^2 f|_{\infty}=\sum_{i,j =1}^{N} |\partial_{x_ix_j}^{2} f|_{L^\infty(\Omega)} $.

\noindent Moreover, we denote

\noindent $\tilde{\Omega}_{r_0}=\Omega \cap B(0, r_0)$,\\
$\Omr=\Omega\setminus (\overline{\omega_0}\cup \overline{\tilde{\Omega}_{r_0}})$,\\
$\Oms=\Omega\setminus \overline{\tilde{\Omega}_{r_0}}$.\\

\subsubsection{The choice of the weight $\sigma$}

In view of the definition of $\psi$ above, we propose the weight \be\label{weight}
    \sigma(t,x) = \theta(t)
    \left(C_\l- |x|^2\psi - \catt\phi\right),
    \quad \phi:=e^{\lambda \psi (x)},
\ee
 where $\lambda$ is a positive parameter aimed to be large and
$r_0$ is a fixed positive constant (small enough) such that it
verifies
\begin{align}\label{condr0}
r_0 &\leq \min \Bigg\{ 1, \frac{\beta_0}{2}, \frac{1}{(2-\gamma)(|D\psi|_{\infty}
+|D^2\psi|_{\infty})},  \frac{1}{\sqrt{2(3|D\psi|_{\infty}^{2}+ |D^2 \psi|_{\infty})}}, \frac{1}{|D\psi|_{\infty}\sqrt{8|\psi|_{\infty}^2+2}},  \nonumber\\
&, \left(\frac{C_3}{8|D\psi|_{\infty}^2
+8|D^2\psi|_{\infty}^2}\right)^{1/(\gamma-1)},
\left(\frac{C_3}{|\mu|
|D\psi|_{\infty}}\right)^{1/(\gamma-1)}, \frac{1}{\sqrt{3|D\psi|_{\infty}}}, \frac{1}{2 |\psi|_\infty |D\psi|_{\infty}}, \nonumber\\
& \frac{1}{\sqrt{8D_{\Omega, \psi_1}|D\psi|_{\infty}/\delta_0 +3 |D^2 \psi|_{\infty}}},   \frac{2}{4 |D\psi|_{\infty} + |D^2 \psi|_{\infty}}
 \Bigg\}.
\end{align}
The normalization by $r_0$ in \eqref{weight} and the election of $r_0$ in \eqref{condr0} are required for
technical reasons needed  later throughout the paper.
 Here $\gamma$   corresponds to the Hardy inequality \eqref{Hardy} with the particular choice $1<  \gamma< 2$,  $C_3$
stands for the constant in inequality \eqref{imprhardy} and  $C_\l$ is
 large enough so as to ensure the positivity of $\sigma$. Besides,
   $\theta$ is defined by
\be\label{theta}
    \theta(t) = \( \frac{1}{t(T-t)}\) ^k,
\ee with $k=1+2/\gamma$.

\subsubsection{Motivation for the choice of  $\sigma$}

Roughly speaking,  the weight $\sigma$ used to prove Carleman estimates for parabolic equations has the general form
 $\sigma(t, x)=\theta(t) A(x)$. In our case, the major difficulty is to match a proper  $A(x)$ because we deal with an equation which differs from the standard heat equation by a singular term in the $x$-variable.
  A positive weight of the form $\sigma_1=\theta(t)(C_\la -e^{\l \psi})$ allows us to control the  heat equation  using a function $\psi$ as in    Fursikov-Imanuvilov \cite{fursikov}.
  Then, this standard weight was modified in Ervedoza \cite{sylvain} when considering  the heat equation with interior quadratic singularity. Basically, the author in \cite{sylvain} proposed a weight which behaves like  $ \sigma_2 \sim \theta(t)(C_\la-|x|^2-|x|^\l)$ as $x$ tends to zero, whereas far away from the origin it still maintains the properties of the standard $\sigma_1$.
The modification near the origin is motivated by some critical terms which must be absorbed outside $\omega$ in the Carleman estimates (see Lemma \ref{I}), i.e.
\begin{multline}\label{iieq}
            + s \is |\partial_n z|^2\  \partial_n \sigma \ds \dt
         +\left(-2 s \iiQ D^2 \sigma(\nabla z, \nabla z) \dx\dt
            + 2 \mu s \iiQ
         \frac{z^2}{|x|^4} x \cdot \nabla \sigma \dx\dt\right)
        \\
        - s\iiQ |\nabla z|^2 \Delta \sigma \ \alpha \dx \dt + \mu s \iiQ \frac{z^2}{|x|^2} \Delta\sigma
         \ \alpha \dx \dt
        - 2 s^3 \iiQ z^2 D^2\sigma \left( \nabla \sigma, \nabla \sigma\right) \dx\dt
         \\
         +s^3\iiQ \alpha z^2 \Delta\sigma
         |\nabla \sigma|^2\ \dx\dt.
     \end{multline}
 In order to take advantage of the optimal Hardy inequality we need to get rid of the singular term  $x \cdot \n \sigma /|x|^4$ in \eqref{iieq} and to impose the   degeneracy $\n \sigma \sim x$ as $x\rightarrow 0$. This fact is reflected in the election of $\sigma_2$ above. However, $\sigma_2$ in \cite{sylvain} does not fit to our case since the move of the singularity  from interior up to the boundary,  produces a loss of regularity for $\sigma_2$, and moreover,  the boundary term in \eqref{iieq} cannot be absorbed in a neighborhood of the origin. For those reasons, we propose the smooth weight $\sigma$ as in \eqref{weight} which makes the terms in \eqref{iieq} positive outside $\omega$,  checking some necessary conditions for $\lambda$ large enough:
\begin{enumerate}
\item $\n \sigma \cdot \vec{n}\geq 0$,  for all $x\in \Gamma$.
\item\label{c2} $-D^2\sigma(x)(\xi, \xi)> 0$,
for all $\xi \in \rr^N$, $x\in \Omega\setminus \overline{B(0, r_0)} $, where $r_0$
verifies \eqref{condr0}.
\item\label{c3} $-\D \sigma > 0 $, for all $x\in \Omega \setminus
\overline{\omega_0}$.
\end{enumerate}

\subsection{Main result} We claim that
\begin{thm}\label{teoremnewvariable}
There exist  positive constants $K$ and $\l_0$ such that for
$\l\geq \l_0$ there exists $s_0=s_0(\l)$ such that for any $s\geq s_0$
we have
\begin{multline}\label{newcarleman}
s \l^2 \iiQr \theta \catt \phi e^{-2s\sigma}|\n w|^2 \dx \dt + s
\iiQ \theta e^{-2s\sigma}\left ( |x|^{2-\gamma} |\n w|^2 +
\frac{w^2}{|x|^\gamma}\right)\dx \dt \\
+ s^3 \iiQ \theta^3 e^{-2s\sigma}|x|^2w^2 \dx\dt + s^3 \l^4 \iiQr
\theta^3e^{-2s\sigma}
\left(\frac{|x|}{r_0}\right)^{3\l}\phi^3 w^2 \dx\dt \\
\leq K \left( s\l^2 \iiq \theta \catt \phi e^{-2s\sigma}|\n
w|^2\dx\dt+ s^3 \l^4 \iiq \theta^3
\left(\frac{|x|}{r_0}\right)^{3\l}\phi^3 e^{-2s\sigma}w^2 \dx\dt
\right).
\end{multline}
\end{thm}

From Theorem  \ref{teoremnewvariable} we can easily obtain the
observability inequality \eqref{obs} via Cacciopoli's inequality.
The details could be reproduced step by step  as in Section 2.2, p.
12, \cite{sylvain}.

\subsection{Preliminaries and useful lemmas}\label{sec:1}

Now, let us assume that $w$ is a solution of \eqref{adjoint} for
some initial data $w_T \in H_0^1(\Omega)$, and define
\be\label{newvar}
    z(x, t) = e^{-s\sigma(x, t)} w(x, t),
\ee which obviously satisfies \be\label{notimeboundary}
    z(x, T) = z(x, 0) = 0, \quad  \mathrm{in}\  H_0^1(\Omega),
\ee due to the assumptions \eqref{hypsigma} on $\sigma$. The
positive parameter $s$ in \eqref{newvar} is meant to be large. Then,
plugging $w = z e^{s\sigma(t,x)}$ in the equation \eqref{adjoint},
we obtain that $z$ satisfies
\begin{multline} \label{eqz}
    \partial_t z + \Delta z + \frac{\mu}{|x|^2} z + 2 s \nabla z\cdot \nabla
    \sigma + s z \Delta
    \sigma
    +z \left( s\partial_t \sigma + s^2|\nabla \sigma|^2 \right) = 0, \quad (x,t)
    \in \Omega\times (0,T),
\end{multline}
with the boundary condition \be\label{nospaceboundary}
            z(x, t) =0, \quad (x, t) \in \Gamma \times (0, T).
\ee

Next, let us define a smooth positive radial function $\alpha(x) =
\alpha(|x|)$ such that \be\label{alpha}
 \alpha(x)=\left\{\begin{array}{ll}
        0, &  |x|\leq r_0/2,\\
      1/N,  & |x|\geq r_0,
         \end{array}\right.
\ee where $r_0>0$ is selected as in  \eqref{condr0}.

Setting \be\label{symanti}
    \left. \begin{array}{ll}
        \displaystyle Sz = \Delta z + \frac{\mu}{|x|^2}
        z + z \left(s \partial_t \sigma + s^2|\nabla \sigma|^2 \right), \\
        Az = \partial_t z + 2 s\nabla z \cdot \nabla
        \sigma + s z\Delta \sigma \left(1+\alpha\right),\\
        Pz=-s \alpha \D \sigma z.
    \end{array}\right.
\ee
one easily deduces from \eqref{eqz} that
$$
    S z + Az +Pz=0, \qquad \quad ||Sz||^2+ ||Az||^2
     + 2 <Sz, Az>  = ||Pz||^2,
$$
where $||\cdot||$ denotes the $L^2((0,T)\times \Omega)$ norm and
$<\cdot,\cdot>$ the corresponding scalar product. In particular, the
quantity \be\label{reducedtothemax}
    I  = <Sz, Az> - \frac{1}{2} ||s \alpha z \Delta\sigma||^2
\ee is non positive.

\begin{lem}\label{I}The following equality holds:
    \begin{multline}\label{Ieq}
        I = -2 s \iiQ D^2 \sigma\left(\nabla z, \nabla z\right) \dx\dt
            + s \is |\partial_n z|^2\  \partial_n \sigma \ds \dt
        \\
        - s\iiQ |\nabla z|^2 \Delta \sigma \ \alpha \dx \dt
        + \frac{s}{2}\iiQ z^2 \Delta^2 \sigma \left(1+\alpha\right) \dx \dt
        \\
        + s\iiQ z^2 \nabla \alpha \cdot \nabla \Delta\sigma \dx \dt
         + \frac{s}{2} \iiQ z^2 \Delta\sigma \ \Delta \alpha \dx \dt
        \\
         - \frac{1}{2} \iiQ z^2 \left( s\partial_{tt}^2 \sigma + 2s^2
         \partial_t \left(|\nabla \sigma|^2 \right)\right)  \dx \dt
     - 2 s^3 \iiQ z^2 D^2\sigma \left( \nabla \sigma, \nabla \sigma\right) \dx\dt
         \\
         +\iiQ \alpha z^2 \Delta\sigma \left(s^2\partial_t \sigma +s^3
         |\nabla \sigma|^2\right) \dx\dt
        -\frac{s^2}{2} \iiQ \alpha^2 z^2 |\Delta \sigma|^2 \dx \dt
        \\
        + \mu s \iiQ \frac{z^2}{|x|^2} \Delta\sigma
         \ \alpha \dx \dt + 2 \mu s \iiQ
         \frac{z^2}{|x|^4} x \cdot \nabla \sigma \dx\dt,
    \end{multline}
    where $\partial_n = \vec{n} \cdot \nabla $ and 
 $\ds$ denotes  the Lebesgue measure on $\Gamma$.
\end{lem}
Here we omit the proof of Lemma \ref{I} since it may be found in
\cite{sylvain}. It is worth mentioning  that the upcoming
computations justified by integrations by parts are done formally.
However, we notice that the final estimates make sense in our
functional framework.  A priori the regularity of the operator
$A_\mu:=-\D - \mu/|x|^2 + \mathcal{C}_{0}^{\gamma} I$ is not enough to justify
the integration by parts since the lack of regularity appears at the
singular point $x=0$.
This issue is presented in a detailed manner in  \cite{cristisumJFA} in the
context of the wave equation with singular potential localized on
the boundary.

Now, we will decompose the term $I$ in \eqref{Ieq} into several
terms that we handle separately.

 Let us define the boundary
term in identity \eqref{Ieq}:

\begin{equation}\label{bdterm}
I_{bd}= \is s |\partial_n z|^2\  \partial_n \sigma \ds \dt.
\end{equation}
 Then define $I_l$  as the sum of the integrals
linear in $\sigma$ which do not have any time derivative:
\begin{multline}\label{il}
    I_l = -2 s \iiQ D^2 \sigma(\nabla z, \nabla z) \dx\dt
        - s \iiQ |\nabla z|^2 \Delta \sigma \ \alpha \dx \dt +
         \frac{s}{2}\iiQ z^2 \Delta^2 \sigma \left(1+\alpha\right) \dx \dt
        \\
        + s \iiQ z^2 \nabla \alpha \cdot \nabla \Delta\sigma \dx \dt +
         \frac{s}{2} \iiQ z^2 \Delta\sigma \ \Delta \alpha \dx \dt
        \\
        + \mu s \iiQ \frac{z^2}{|x|^2} \Delta \sigma  \
        \alpha \dx \dt + 2 \mu s \iiQ \frac{z^2}{|x|^4} x\cdot
        \nabla\sigma \dx \dt.
\end{multline}
We then consider the sum of the integrals involving non-linear terms in
$\sigma$ and without any time derivative, that is
\begin{multline}\label{Inl}
I_{nl} =  - 2 s^3 \iiQ z^2 D^2\sigma \left( \nabla \sigma, \nabla
\sigma\right) \dx\dt
         + s^3 \iiQ \alpha z^2 \Delta\sigma |\nabla \sigma|^2 \dx\dt \\
        -\frac{s^2}{2} \iiQ \alpha^2 z^2 |\Delta \sigma|^2\dx.
\end{multline}
We finally estimate the terms involving the time derivatives in
$\sigma$:
\begin{align}\label{it}
    I_t =  - \frac{1}{2} \iiQ z^2 \left( s \partial_{tt}^2 \sigma + 2 s^2 \partial_t \left(|\nabla \sigma|^2 \right)\right)
\dx \dt
         + s^2 \iiQ \alpha z^2 \Delta\sigma \partial_t \sigma\dx\dt.
\end{align}
In the next step we give convenient estimations for the
terms defined above. In order to do that several lemmas are proved.

\begin{lem}\label{pozbd}
It holds that $I_{bd}\geq 0$, for any $\l>0$.
\end{lem}

\begin{lem}\label{lemaIll}
There exists $\la_0$ such that for any $\la\geq \la_0$ and any $s>0$
then
\begin{multline}
I_l \geq C_{3} s  \iiQ \theta \left( |x|^{2-\gamma}|\n z|^2 +
\frac{z^2}{|x|^\gamma}\right)\dx \dt +C_8 s\l^2 \iiQr \theta \catt \phi |\n z|^2 \dx \dt+ \\
+ s\la \iir0 \theta \left(\frac{|x|}{r_0}\right)^{\l-2}|\n z|^2 \dx
\dt 
- B_\la s \iiQ \theta z^2 \dx\dt\\
 - C_7 s\la^2 \iiq \theta \catt
\phi |\n z|^2 \dx \dt,
\end{multline}
where the constants $C_3, C_8, C_7$  are uniform  in $s$
and  $\la$, and $B_\l$ is uniform in $s$.
\end{lem}

\begin{lem}\label{Inlfinal}
There exists $\la_0$ such that for any $\la\geq \la_0$ there exists
$s_0=s_0(\la)$ such that for any $s\geq  s_0$ it holds
\begin{multline}
 I_{nl}\geq \frac{s^3}{2} \iir0 \theta^3 |x|^2 z^2
\dx \dt+ C_{15} s^3 \l^4 \iiQr \theta^3
\left(\frac{|x|}{r_0}\right)^{3\l} \phi^3
z^2 \dx \dt\\
- C_{16} s^3\l^4 \iiq \theta^3 \left(\frac{|x|}{r_0}\right)^{3\l}\phi^3
 z^2 \dx\dt,
\end{multline}
for some constants $C_{15}$, $C_{16}$ uniform in $s$ and $\la$.
\end{lem}

Taking into account the negative terms in the expression of $I_l$
that we want to get rid of, we define
\begin{equation}
I_r=I_t
-B_\la s \iiQ \theta z^2 \dx\dt.
\end{equation}
Then
\begin{lem}\label{Ir}
There exists $\la_0$ such that for any $\la\geq \la_0$ there exists
$s_0=s_0(\la)$ such that for any $s\geq  s_0$ we have
\begin{multline}
|I_r|\leq \frac{C_{3} s}{2} \iiQ \theta \frac{z^2}{|x|^\gamma} \dx
\dt + \frac{C_{15}}{2}s^3 \l^4 \iiQs \theta^3
\left(\frac{|x|}{r_0}\right)^{3\l} \phi^3  z^2 \dx \dt +
\frac{s^3}{4} \iir0 \theta^3 |x|^2 z^2 \dx \dt.
\end{multline}
\end{lem}

From the lemmas above we obtain the Carleman inequality in the  variable $z$ as follows.
\begin{thm}\label{teoremnewvariable1}
There exist positive constant $K$ and $\l_0$ such that for
$\l\geq \l_0$ there exists $s_0=s_0(\l)$ such that for any $s\geq s_0$
we have
\begin{multline}\label{newcarleman}
s \l^2 \iiQr \theta \catt \phi |\n z|^2 \dx \dt + s  \iiQ \theta
\left( |x|^{2-\gamma}|\n z|^2 +
\frac{z^2}{|x|^\gamma}\right)\dx \dt \\
+ s^3 \iiQ \theta^3 |x|^2z^2 \dx\dt + s^3 \l^4 \iiQr \theta^3
\left(\frac{|x|}{r_0}\right)^{3\l}\phi^3 z^2 \dx\dt \\
\leq K \left( s\l^2 \iiq \theta \catt \phi |\n z|^2\dx\dt+ s^3 \l^4
\iiq \theta^3 \left(\frac{|x|}{r_0}\right)^{3\l} \phi^3 z^2 \dx\dt \right).
\end{multline}
\end{thm}

Coming back from the variable  $z$ to the  solution $w$, due to \eqref{newcarleman} we obtain the
conclusion of Theorem \ref{teoremnewvariable}.

\section{Basic computations}\label{babizduca}
This section is based on some preliminary computations which will be applied in Section \ref{sec:1} to the proofs of Lemmas \ref{lemaIll}, \ref{Inlfinal}, \ref{Ir}.
\subsection{Fundamental property of $\psi$.}
The main result of this paragraph is
\begin{lem}\label{crucial}
Assume $\psi$ is the weight defined in \eqref{defimp} by means of $\psi_1$ and $\delta$. Then  there exists a constant $D_{\Omega, \psi_1} >0$,  which depends on $\Omega$ and $\psi_1$,  such that
\begin{equation}\label{veryimp}
|x \cdot \n \psi (x) - \delta \psi_1(x)| \leq \delta D_{\Omega, \psi_1} |x|^2, \quad \forall x\in \Omega.
\end{equation}
\end{lem}
This lemma consists in a basic result which is worth mentioning  since it plays a crucial role in the proof of Theorem \ref{teoremnewvariable1}. Of course, Lemma \ref{crucial} makes sense to be proven close to the origin, otherwise it is a triviality.
\begin{proof}[Proof of Lemma \eqref{crucial}]
We split the proof in two important steps.\\

\noindent \textit{Step 1:}  There exists a constant $E_\Omega>0$ such that
\be\label{vrimp2}
|\pr | \leq E_\Omega |x|, \quad \textrm{for } x\in \Omega_{\beta_0} \textrm{  near the origin}.
\ee
Indeed, without loss of generality we may assume that the parametrization of $\Gamma$ near the origin is given by $x_N=\beta|x'|^2$, for some $\beta\in \rr$.
Then, for a fixed point  $x\in \Omega_{\beta_0}$ near the origin, its projection on $\Gamma$ in \eqref{puffy} is given by $\pr=(a', \beta |a'|^2)$ and minimizes the functional
\begin{equation}
f(x; a')=||x'-a'||^2+||x_N-\beta |a'|^2||^2.
\end{equation}
In other words,  $\pr$ verifies $\n_{a'} f=0$ which is equivalent to
\begin{equation}\label{idd}
2(a'-x')-4\beta a'(x_N-\beta |a'|^2)=0.
\end{equation}
Multiplying \eqref{idd} by $a'\not \equiv 0$  we have
$$2|a'|^2-2<x', a'> -4\beta |a'|^2 x_N +4\beta^2 |a'|^4=0,$$
and therefore from Cauchy-Schwartz inequality we get
\begin{align}\label{lin1}
4\beta^2 |a'|^4 &= -2|a'|^2 +2 <x', a'> +4\beta |a'|^2 x_N \nonumber\\
&\leq 2 |x'||a'| +4|\beta| |a'|^2 |x_N|\nonumber\\
& \leq \max\{ 2, 4|\beta| R_\Omega\} |a'|(|x'|+|x_N|)\nonumber\\
 &\leq \sqrt{2} \max\{ 2, 4|\beta| R_\Omega\} |a'||x|.
\end{align}
From \eqref{lin1} we deduce
\begin{equation}\label{lin2}
4\beta^2 |a'|^3 \leq \sqrt{2} \max\{ 2, 4|\beta| R_\Omega\} |x|.
\end{equation}
On the other hand, from \eqref{idd} we have
\begin{align}\label{nice1}
2|a'|&\leq 2|x'| +4|\beta| |a'||x_N| +4\beta^2 |a'|^3\nonumber\\
&\leq \max\{ 2, 4|\beta| R_\Omega\} (|x'|+|x_N|) + 4\beta^2 |a'|^3\nonumber\\
&\leq \sqrt{2}\max\{ 2, 4|\beta| R_\Omega\} |x|+4\beta^2 |a'|^3,
\end{align}
which combined with \eqref{lin2} leads to
\begin{equation}\label{lin3}
|a'|\leq \sqrt{2}\max\{2, 4|\beta| R_\Omega\} |x|.
\end{equation}
Formula \eqref{lin3} yields to
\begin{align}
|\pr|^2&=|a'|^2+\beta^2 |a'|^4\nonumber\\
& \leq (1+ \beta^2 R_\Omega^2)|a'|^2\nonumber\\
&\leq 2(1+\beta^2R_\Omega^2) \max\{ 2, 4|\beta|R_\Omega\}^2 |x|^2,
\end{align}
which concludes the proof of \eqref{vrimp2}. \\

\noindent \textit{Step 2:}
According to the definition \eqref{defimp} of $\psi$ and the properties of the distance function $\rho(x)$, for $x\in \Omega_{\beta_0}$ we have
 \begin{align}\label{pru}
 |x\cdot \n \psi(x) -\delta \psi_1(x) |&= \delta |x\cdot \n \rho(x)-\rho(x)|\nonumber\\
 &= \delta |(\pr - \rho(x)\vec{n}(\pr))(-\vec{n}(\pr))-\rho(x)|\nonumber\\
 &= \delta |\pr \cdot \vec{n}(\pr)|.
 \end{align}
In addition,   applying \eqref{pru}, \eqref{impcond1} and \eqref{vrimp2} we get
 \be\label{impcond}
 |x\cdot \n \psi(x) -\delta \psi_1(x) | \leq \delta C_\Omega|\pr|^2 \leq \delta C_\Omega E_\Omega^2  |x|^2, \quad \forall x\in \Omega_{\beta_0} \cap B(0, \nu_0),
 \ee
 for some $\nu_0>0$ small enough. Due to \eqref{impcond}
 the proof of Lemma \ref{crucial} is obtained in a neighborhood of the origin for $D_{\Omega, \psi_1}=C_\Omega E_\Omega^2$ (close to the origin, the dependence on  $\psi_1$ is involved in the identification $\psi_1=\rho$). Far way from the origin, the proof of \eqref{impcond} is trivial with $D_{\Omega, \psi_1}$ depending on $\Omega$, $|\psi_1|_{\infty}$ and  $|D\psi_1|_{\infty}$. Thus, the proof of Lemma \ref{crucial} is finished.
\end{proof}
\subsection{Useful identities}

 Part of the
computations here require a more careful analysis. First of all, for $\sigma$ as in \eqref{weight},
 we make the notations
$$\sigma_{x^2}=-\theta(t) \tau_{x^2}, \textrm{ where } \tau_{x^2}=|x|^2\psi,$$
respectively
$$\sigma_\phi=-\theta(t)  \tau_\phi, \textrm{ where }
\tau_\phi=
\catt\phi.
$$
and $$\tau =\tau_{x^2}+\tau_\phi.$$
Next we deduce some formulas for $\tau_{x^2}$ and $\tau_{\phi}$ that
we are going to use in our computations. More precisely, for all $x
\in \rr^N$ and any $i, j \in \{1, \ldots, N\}$
  we have
\begin{align}\label{catevaformule}
\p_{x_i} \tau_{x^2}&=2x_i\psi + |x|^2 \p_{x_i} \psi,\nonumber\\
\partial_{x_i x_j}^{2}\tau_{x^2}&=2\delta_{i,j}\psi + 2 x_i \partial_{x_j}\psi+2
x_j \partial_{x_i}\psi +|x|^2 \partial_{x_i x_j}^{2} \psi,
\end{align}
and
\begin{align}
\D\tau_{x^2}&= 2N \psi + 4(x\cdot \n \psi) + |x|^2 \D \psi, \label{id1:1}\\
D^2\tau_{x^2}(\xi, \xi)&= 2\psi |\xi|^2+ 4(x\cdot \xi)(\n \psi \cdot \xi)+ |x|^2 D^2 \psi (\xi, \xi). \label{id1:2}
\end{align}
On the other hand,
\begin{align}\label{form2}
  \p_{x_i} \tau_\phi &=\caat\big(\lambda x_i |x|^{\lambda-2}
  +\lambda  |x|^{\lambda}\p_{x_i} \psi
\big)\phi,\nonumber\\
 \p^2_{x_ix_j}\tau_\phi&=  \caat\Big(\lambda
\delta_{ij}|x|^{\l-2}+\l(\lambda-2) x_i x_j |x|^{\lambda-4}+\lambda^2 x_j
\p_{x_i}\psi |x|^{\l-2}+\l^2
x_i \p_{x_j} \psi |x|^{\l-2} \nonumber\\
&+\l |x|^{\l} \p_{x_ix_j}\psi +\l^2 |x|^{\l} \p_{x_i}\psi
\p_{x_j}\psi \Big)\phi,
\end{align}
and in particular
\begin{align}
\D \tau_{\phi} &= \caat \Big( \l(\l+N-2)|x|^{\l-2}+2\l^2 (x\cdot \n \psi)|x|^{\l-2} + \l \D \psi |x|^\l  + \l^2 |\n \psi|^2 |x|^\l \Big) \phi, \label{eeq:1}\\
D^2 \tau_{\phi }(\xi, \xi) &=\caat \Big(\l |\xi|^2 |x|^{\l-2} + \l(\l-2)|x\cdot \xi|^2 |x|^{\l-4} + 2\l^2 (x\cdot \xi)(\n \psi \cdot \xi)|x|^{\l-2}\nonumber\\
 &+ \l D^2\psi (\xi,\xi)|x|^\l + \l^2 |\n \psi \cdot \xi |^2 |x|^{\l}\Big )\phi. \label{eeq:2}
\end{align}

\subsection{Upper and lower bounds for $\D \tau_{x^2}$, $\D \tau_{\phi}$, $D^2 \tau_{x^2}(\xi,\xi)$, $D^2 \tau_{\phi}(\xi,\xi)$}
\begin{prop}\label{prop3}
For $r_0$ as in \eqref{condr0} we have
\begin{align}
\Delta \tau_{x^2} \geq 0, \  D^2\tau_{x^2}(\xi, \xi) &\geq 0, \qquad \forall x\in \tilde{\Omega}_{r_0}, \ \forall \xi \in \rr^N,  \label{eq:1}\\
| D^2\tau_{x^2}(\xi, \xi)| &\leq D_1 |\xi|^2, \qquad \forall x\in \Omega, \ \forall \xi \in \rr^N,\label{eq:2}\\
 |\Delta \tau_{x^2}| &\leq D_1, \qquad \forall x\in \Omega, \ \forall \xi \in \rr^N, \label{eq:3}
\end{align}
where $D_1$ is a large enough constant depending on $\Omega$ and $\psi$.
\end{prop}
\begin{prop}\label{prop4}
For $r_0$ as in \eqref{condr0} and any $\lambda\geq 6$ we have that
\begin{align}
D^2 \tau_\phi (\xi, \xi) &\geq  \frac{\lambda}{2} \left (\frac{|x|}{r_0}\right)^{\l-2} \phi |\xi|^2, \qquad \forall x\in \tilde{\Omega}_{r_0}, \ \forall \xi \in \rr^N, \label{eqq:1}\\
\D \tau_{\phi} &\geq \l^2 \catt  \phi,  \qquad \forall x\in \Omega\setminus \overline{\omega_0}, \label{eqq:3}\\
D^2 \tau_{\phi}(\xi,\xi) &\geq -\l D_4 \ccat \phi, \qquad \forall x\in \Omega,\label{eqq:4}
\end{align}
where the constant $D_4=D_4(\Omega, r_0, |D\psi|_{\infty})$ is  uniform in $\la$.
\end{prop}
\begin{proof}[Proof of Proposition \ref{prop3}]
Observe that the proofs of \eqref{eq:2} and \eqref{eq:3} are consequences of the $C^2$ regularity of $\tau_{x^2}$.
To conclude, it is enough to show that $D^2\tau_{x^2}(\xi, \xi)\geq 0$ since this also implies that $\D \tau_{x^2}\geq 0$ (we just have to choose $\xi=e_i$, for all  $i\in \{1, \ldots, N\}$).
Indeed, from \eqref{id1:2} and Cauchy-Schwartz inequality we obtain
\begin{align}\label{ineq1}
D^2\tau_{x^2}(\xi, \xi) & \geq |\xi|^2 (2\psi- 4|x| |\n \psi|- |x|^2 |D^2 \psi |_{\infty}).
\end{align}
Since $\psi(x)\geq 1$ for $x\in \tilde{\Omega}_{r_0}$, from \eqref{ineq1} we finally get
\begin{align}
D^2\tau_{x^2}(\xi, \xi)
&\geq |\xi|^2 (2 -4 r_0 |D\psi|_{\infty} - r_0^2 |D^2 \psi|_{\infty})\nonumber\\
& \geq |\xi|^2 \big(2 - r_0 (4 |D\psi|_{\infty} +  |D^2 \psi|_{\infty})\big)\nonumber\\
& \geq 0,
\end{align}
 since $r_0$ satisfies \eqref{condr0}.
\end{proof}
\begin{proof}[Proof of Proposition \ref{prop4}]
Firstly,  in \eqref{eeq:2} we write $D^2 \tau_\phi(\xi, \xi)=\phi (1/r_0)^{\la}
S_\phi$ where
\begin{align}\label{form4}
S_\phi&=\l |\xi|^2 |x|^{\l-2}+ \l (\l-2) |x\cdot \xi|^2 |x|^{\l
-4}+2\l^2 (x\cdot \xi)(\xi \cdot \n
\psi)|x|^{\l-2}\nonumber\\
&+\l |x|^{\l} D^2\psi (\xi, \xi) +\l^2 |x|^{\l }|\n \psi \cdot \xi
|^2.
\end{align}
Next, we have the inequality
$$\Big|2\l^2 (x\cdot \xi)(\n \psi \cdot \xi)|x|^{\l-2}\Big|\leq a \l^2|x\cdot \xi|^2 |x|^{\l-4}+
\frac{\l^2}{a}|x|^{\l}|\n \psi \cdot \xi|^2, \q \forall  a>0, $$
which combined with \eqref{form4}  leads to
\begin{align*}
S_\phi\geq \l |\xi|^2|x|^{\l-2} + (\l^2 -2\l -a \l^2)|x\cdot \xi|^2
|x|^{\l-4}+ \l |x|^{\l}D^2\psi(\xi, \xi)+ \left(\l^2
-\frac{\l^2}{a}\right)|x|^{\l}|\n \psi\cdot \xi |^2.
\end{align*}
Choosing $a>0$ such that $\l^2-2\l -a \l^2=0$ (i.e.
$a=(\l-2)/\l$), we remark that
\begin{align}\label{form5}
S_\phi&\geq \l |\xi|^2 |x|^{\l-2}+ \l |x|^{\l}D^2\psi (\xi,\xi)-
\frac{2\l^2}{\l-2}|x|^{\l}|\n \psi|^2 |\xi|^2, \q \forall x, \ \forall
\xi.
\end{align}
Applying \eqref{form5} for $x\in \tilde{\Omega}_{r_0}$ and $\l\geq 6$,  we deduce
\begin{align}
S_\phi & \geq \frac{\l}{2}|x|^{\l-2}|\xi|^2 + \l|x|^{\l-2}|\xi|^2
\left(\frac{1}{2}-\frac{2
\l }{\l-2}|x|^2|D\psi|_{\infty}^{2}-|x|^2|D^2\psi|_{\infty}\right)\nonumber\\
& \geq \frac{\l}{2}|x|^{\l-2}|\xi|^2 + \l|x|^{\l-2}|\xi|^2
\left(\frac{1}{2}-r_0^2\left(
3|D\psi|_{\infty}^{2}+|D^2\psi|_{\infty}\right)\right)\nonumber\\
&\geq  \frac{\l}{2}|x|^{\l-2}|\xi|^2,
\end{align}
which holds true for $r_0$ as in \eqref{condr0}.
 This yields the
proof of \eqref{eqq:1}.

Next, let us prove \eqref{eqq:3}.

\noindent According to Lemma \ref{veryimp}, the definition of $\psi$ and  \eqref{eeq:1} we get
  \begin{align}
\D \tau_{\phi}&\geq \frac{1}{r_0^\l} \Big( \l(\l+N-2)|x|^{\l-2} + 2 \l^2 (\delta\psi_1 -\delta D_{\Omega, \psi_1} |x|^2)|x|^{\l-2}+ \l \D \psi |x|^\l + \l^2 |\n \psi|^2 |x|^\l \Big)\phi\nonumber\\
& \geq \frac{1}{r_0^\l} \Big( \l(\l+N-2)|x|^{\l-2}- 2 \l^2 \delta D_{\Omega, \psi_1} |x|^{\l} +  \l \D \psi |x|^\l + \l^2 |\n \psi|^2 |x|^\l \Big)\phi\nonumber\\
& \geq \l^2 \catt \phi \left(|\n \psi|^2- 2 \delta D_{\Omega, \psi_1} -\frac{|\D \psi|}{\lambda}\right)\nonumber\\
&\geq \l^2 \catt \phi \left(\delta^2\delta_0^2- 2\delta D_{\Omega, \psi_1} -\delta |D^2\psi_1|_{\infty}\right)\nonumber\\
& \geq \l^2 \catt \phi, \quad \forall x\in \Omega\setminus \overline{\omega_0},
\end{align}
provided $\delta$ is large enough $(\delta \geq \max\{(1+2D_{\Omega, \psi_1}+ |D^2\psi_1|_{\infty})/\delta_0^2, 1\})$ and $\l\geq 1$.

For the proof of \eqref{eqq:4} we proceed as follows.
We observe that
$$\l^2 |x\cdot \xi|^2 |x|^{\l-4}+ 2\l^2 (x\cdot \xi)(\n  \psi\cdot \xi)|x|^{\l-2}+ \l^2 |\n \psi\cdot \xi|^2 |x|^\l \geq \l^2 \left(|x\cdot \xi||x|^{\l/2-2}-|\n \psi \cdot \xi||x|^{\l/2} \right)^2 \geq 0.$$
 This with \eqref{eeq:2} give
\begin{align}
D^2 \tau_{\phi} (\xi, \xi) & \geq  \frac{1}{r_0^\l}\left(-3\l |\xi|^2 |x|^{\l-2}+\l D^2\psi(\xi, \xi)|x|^\l \right) \phi\nonumber\\
&\geq -\l |\xi|^2\ccat \left(\frac{3}{r_0^2} +  |D^2\psi|_{\infty} \frac{|x|^2}{r_0^2} \right)\phi, \quad \forall x\in \Omega,
\end{align}
which concludes the validity of \eqref{eqq:4} for $D_4=(3+|D^2\psi_1|_{\infty}R_\Omega^2)/r_0^2$.
\end{proof}

\subsection{Bounds for $2D^2\tau (\n \tau, \n \tau )- \alpha \D \tau |\n \tau|^2$}

In this subsection we provide very useful pointwise estimates for the term
$$2D^2\tau (\n \tau, \n \tau )- \alpha \D \tau |\n \tau|^2,$$
 which appears in identity \eqref{Ieq} of Lemma \ref{I}.  These computations represent the most technical part of the paper. Besides that, they turn out to play a crucial role in proving  Carleman estimates and observability.
Before going into details we need some a priori technical identities.

Firstly, due to \eqref{catevaformule}-\eqref{form2} we have
\begin{align}
\p_{x_i}\tau&= 2x_i\psi + |x|^2 \p_{x_i} \psi+ \bou\Big(\lambda x_i
|x|^{\lambda-2}
  +\lambda  |x|^{\lambda}\p_{x_i} \psi
\Big), \label{prost1}\\
\p^2_{x_ix_j}\tau&=2\delta_{ij}\psi + 2 x_i \p_{x_j}\psi+2 x_j
\p_{x_i}\psi +|x|^2 \p_{x_i x_j}^{2} \psi\nonumber\\
&+\bou\Big(\lambda
\delta_{ij}|x|^{\l-2}+\l(\lambda-2) x_i x_j
|x|^{\lambda-4}+\lambda^2 x_j \p_{x_i}\psi |x|^{\l-2}\nonumber\\
&+\l^2 x_i \p_{x_j} \psi
|x|^{\l-2}+\l |x|^{\l} \p_{x_ix_j}\psi +\l^2 |x|^{\l} \p_{x_i}\psi
\p_{x_j}\psi \Big), \label{otherform}
\end{align}

and in consequence
\begin{align}
\D \tau&=2 N \psi + 4 (x \cdot \n \psi)|x|^2 \D  \psi\nonumber\\
&+\bou \Big(\l(N+\lambda-2)
|x|^{\lambda-2}+2\lambda^2 (x\cdot \n \psi) |x|^{\l-2} +\l |x|^{\l} \D \psi +\l^2 |x|^{\l} |\n \psi
|^2  \Big), \label{tavico1}\\
D^2\tau (\xi, \xi) &=2 \psi|\xi|^2  + 4 (x\cdot \xi)(\n \psi  \cdot \xi)  +|x|^2 D^{2}\psi(\xi, \xi)\nonumber\\
&+\bou \Big(\lambda |x|^{\l-2}|\xi|^2+\l(\lambda-2) |x\cdot \xi |^2
|x|^{\lambda-4}\nonumber\\
&+ 2\lambda^2 (x \cdot \xi) (\n \psi\cdot \xi) |x|^{\l-2}+\l |x|^{\l} D^2\psi (\xi, \xi) +\l^2 |x|^{\l} |\n \psi \cdot \xi
|^2 \Big). \label{tavico2}
\end{align}

Using the expressions in \eqref{prost1}-\eqref{otherform} we obtain several useful
formulas:
\begin{align}
|x\cdot \n \tau|^2&= |x|^2|\n \tau|^2+\Big(|\n \psi
\cdot x|^2-|x|^2|\n \psi|^2\Big)\left(|x|^2 +\l \bou |x|^\l\right)^2, \label{frig1} \\
 (x \cdot \n \tau)(\n \psi \cdot \n \tau)&=
 |\n \tau|^2 \n \psi \cdot x+\Big(|x|^2|\n
\psi|^2-|x\cdot \n \psi|^2 \Big) \left(1+ \l \bou
|x|^{\l-2}\right)\times \nonumber\\
&\times\left(2\psi |x|^2 + \l \bou |x|^\l\right), \label{frig2}\\
|\n \psi \cdot \n \tau|^2&=|\n \psi|^2|\n \tau|^2+ \Big(|x\cdot \n
\psi|^2-|x|^2|\n \psi|^2 \Big) \left(2\psi+\l \bou
|x|^{\l-2}\right)^2.\label{frig3}
\end{align}
Using the identities  \eqref{prost1} and   \eqref{frig1}-\eqref{frig3} we conclude
\begin{align}
2D^2\tau (\n \tau, \n \tau )- \alpha \D \tau |\n \tau|^2:=T_1+T_2+T_3,\label{cost1}
\end{align}
where
\begin{align}
T_1&=2 \psi (2-\alpha N)|\n \tau|^2+ 4(2-\alpha)|\n \tau|^2 \n \psi
\cdot x+
2 |x|^2 D^2\psi (\n \tau, \n \tau)-\alpha |x|^2 \D \psi |\n \tau|^2,\label{cost2}
\end{align}
\begin{align}
T_2&=8 \big(|x|^2 |\n \psi|^2-|\n \psi \cdot x|^2\big) \left(1+ \l  \bou
|x|^{\l-2}\right)\left(2\psi|x|^2 +\l \bou  |x|^\l\right)\nonumber\\
&+\frac{\phi}{r_0^\l} \big(|x|^2 |\n
\psi|^2-|\n \psi \cdot x|^2\big) \Bigg( 4\l^3 \left(\frac{1}{r_0}\right)^{2\l} \phi^2 |x|^{3\l -4}
+8 \l^2 \bou |x|^{2\l-2}  \nonumber\\
&+ 8 \l^2 \psi (1-\psi) |x|^\l-
2\l (\l-2) |x|^\l \Bigg),\label{cost3}
\end{align}
\begin{align}
T_3&=\bou \bigg\{ \Big[\left((2-\alpha) \l^2 - \l (2+ \alpha N-2
\alpha ) \right)|x|^{\l-2}+2 \l^2 (2-\alpha)|x|^{\l-2}\n \psi
\cdot x     \nonumber\\
&-\alpha \l |x|^\l \D \psi  +(2-\alpha) \l^2 |x|^{\l} |\n \psi|^2 \Big] |\n \tau|^2+2 \l |x|^\l D^2 \psi (\n \tau, \n \tau)  \bigg\}.\label{cost4}
\end{align}
Based on all this we can claim
\begin{prop}\label{prop8}
For $r_0$ as in \eqref{condr0}, there exist constants $D_5, D_6>0$ depending on $\psi$ and $\Omega$, such that  $T_1$ in \eqref{cost2} satisfies the following bounds:
\begin{align}
T_1 & \geq |\n \tau|^2, \quad \forall x\in \tilde{\Omega}_{r_0},\label{chici1}\\
T_1 & \geq - D_5 |\n \tau|^2, \quad \forall x\in \Omr\label{chici2},\\
|T_1| & \leq D_6 |\n \tau|^2, \quad \forall x\in \omega_0. \label{chici3}
\end{align}

\end{prop}
\begin{prop}\label{prop5}
 There exists $\l_0=\lambda_0(R_\Omega)>0$  such that for any $\l\geq \l_0$ and  $r_0$ as in \eqref{condr0}, the term $T_2$ in \eqref{cost3} verifies
\begin{align}
T_2 & \geq -\bou |D \psi|_{\infty}^{2} (8 \l^2 \psi^2 +2 \l^2 )|x|^{\l+2}, \quad \forall x\in \tilde{\Omega}_{r_0}, \label{bambu1}\\
T_2 & \geq 0, \quad \forall x\in \Oms. \label{bambu2}
\end{align}
\end{prop}
\begin{prop}\label{prop6}
There exists $\l_0$ large enough such that for any $\l\geq \l_0$,  and  $r_0$ as in \eqref{condr0}, 
 the term $T_3$ in \eqref{cost4} verifies
\begin{align}
T_3 & \geq \l^2 \left( \bou |x|^{\l-2} +  \catt \phi \right) |\n \tau|^2, \quad \forall x\in \Omega\setminus \overline{\omega_0},\label{tumpi1}\\
T_3 & \leq  D_7 \l^2 \bou |x|^{\l-2} |\n \tau|^2, \quad \forall x\in \Omega, \label{tumpi3}
\end{align}
for some constant $D_7=D_7(\Omega, \psi)$ uniform in $\l$.
\end{prop}
\begin{prop}\label{tzuca}
For any $\lambda>1$  and  $r_0$ as in \eqref{condr0}, it holds that
\begin{align}
 |\n \tau|^2 &\geq |x|^2, \quad  \forall  x\in \tilde{\Omega}_{r_0},\label{pici}\\
 |\n \tau|^2 &\geq  \l^2 \left(\frac{|x|}{r_0}\right)^{2\l}\phi^2, \quad \forall x\in \Omr, \label{pici2}\\
 |\n \tau|^2 &\leq D_8 \l^2 \left(\frac{|x|}{r_0}\right)^{2\l}\phi^2, \quad \forall x\in \omega_0,   \label{pici3}
\end{align}
where $D_8$ is a constant depending only on $\Omega$ and $\psi$.
\end{prop}
\begin{proof}[Proof of Proposition \ref{prop8}]
Let us firstly prove \eqref{chici1}.  Since $\alpha$ satisfies \eqref{alpha}, due to the properties of $\psi$ and Lemma \ref{veryimp} we obtain
\begin{align*}
T_1 & \geq\Big(2 + 4(2-\alpha) (\delta \psi_1- \delta D_{\Omega, \psi_1} |x|^2) - 2 |x|^2 |D^2 \psi|_{\infty}- |x|^2 |D^2 \psi|_{\infty} \Big) |\n \tau|^2\nonumber\\
& \geq  \Big(2-r_0^2\big(8\delta D_{\Omega, \psi_1} + 3 |D^2 \psi|_{\infty}\big)\Big) |\n \tau|^2 \nonumber\\
& \geq  \left(2-r_0^2\left(8\frac{D_{\Omega, \psi_1}}{\delta_0} | D \psi|_{\infty} + 3 |D^2 \psi|_{\infty}\right)\right) |\n \tau|^2,\nonumber\\
& \geq |\n \tau|^2,   \quad \forall x\in \tilde{\Omega}_{r_0},
\end{align*}
which is true since $r_0\leq 1/(\sqrt{8D_{\Omega, \psi_1}|D\psi|_{\infty}/\delta_0 +3 |D^2 \psi|_{\infty}})$.
On the other hand, inequalities \eqref{chici2}-\eqref{chici3} are obvious due to the $C^2$ regularity of $T_1$.
\end{proof}
\begin{proof}[Proof of Proposition \ref{prop5}]
The proof of \eqref{bambu1} is a consequence of the Cauchy-Schwartz inequality. Besides, \eqref{bambu2} holds true for $\l$ large enough since the term containing $\l^3$ is positive and dominates all the other terms far away from the origin.
\end{proof}
\begin{proof}[Proof of Proposition \ref{prop6}]

Due to the definition of $\alpha$ in \eqref{alpha}, for $\lambda\geq \l_0$ large enough we have
\be\label{incep}
(2-\alpha)\l^2 -\l(2+\alpha N-2 \alpha)\geq \l^2, \quad \forall x\in \Omega.
\ee
 Therefore, from Lemma \ref{veryimp} and the properties of $\psi$, for $x\in \Omega\setminus \overline{\omega_0}$ we have
\begin{align}
T_3&\geq \bou \Big(\l^2 |x|^{\l-2} + 2\l^2 |x|^{\l-2} (\delta \psi_1-\delta D_{\Omega, \psi_1} |x|^2)-\nonumber\\
&-\l |D\psi|_{\infty} |x|^\l + \l^2 |x|^\l |\n \psi|^2 -2\l |D^2 \psi|_{\infty} |x|^\l \Big)|\n \tau|^2\nonumber\\
& \geq \l^2 \bou |x|^{\l-2} |\n \tau|^2 + \bou \Big(\l^2 |x|^\l \delta^2 |\n \psi_1|^{2}\nonumber\\
&-2\delta D_{\Omega, \psi_1} \l^2 |x|^\l - \l \delta |D\psi_1|_{\infty} |x|^\l - 2 \l \delta |D^2 \psi_1|_{\infty} |x|^\l \Big) |\n \tau|^2 \nonumber\\
& \geq \l^2 \bou |x|^{\l-2} |\n \tau|^2+ \l^2 \Big(\delta^2 \delta_0^2 - 2 \delta D_{\Omega, \psi_1} \nonumber\\
&-\frac{\delta}{\l} |D\psi_1|_{\infty} -\frac{2 \delta}{\l} |D^2 \psi_1|_{\infty}\Big) \bou |x|^\l |\n \tau|^2\nonumber\\
& \geq \l^2 \bou |x|^{\l-2} |\n \tau|^2+ \l^2 \bou |x|^\l |\n \tau|^2,\nonumber
\end{align}
  due to $\delta \geq \max\{ (1+ 2 D_{\Omega, \psi_1} + |D\psi_1|_{\infty}+ 2|D^2 \psi_1|_{\infty})/\delta_0^2, 1\}$  and $\la\geq \max\{\l_0, 1\}$.

Again, the proof of \eqref{tumpi3} is a consequence of the $C^2$ regularity of $T_3$.
\end{proof}
\begin{proof}[Proof of Proposition \ref{tzuca}]
Expanding the square in \eqref{prost1} we obtain
\begin{align}\label{gic2}
|\n \tau|^2&=4|x|^2\psi^2+ |x|^4|\n \psi|^2 + \l^2 \left(\bou\right)^2 \Big| x|x|^{\l-2}+|x|^\l \n \psi \Big|^2\nonumber\\
&+ 4\psi |x|^2 (x\cdot \n \psi) + 4\l \psi \bou (|x|^{\l-1}+|x|^\l \n \psi \cdot x)\nonumber\\
&+ 2\l\bou |x|^2 (x\cdot \n \psi|x|^{\l-2} +|x|^\l |\n \psi|^2).
\end{align}
In order to absorb the cross term $x\cdot \nabla \psi$ near the origin, we observe that $$3|x|^2+4\psi |x|^2(\n \psi \cdot x)>0, \quad  4\psi (|x|^{\l-1}+|x|^{\l} \n \psi \cdot x)+2 x\cdot \n \psi |x|^\l>0, \qquad x \in \tilde{\Omega}_{r_0},$$ since $r_0 < 1/ (2 |\psi|_{\infty} |D\psi|_{\infty})$ respectively $r_0 \leq 1/\sqrt{3|D\psi|_{\infty}}$. Hence, we conclude the validity of \eqref{pici}.
For the proof of \eqref{pici2} we proceed in several steps. First, let us observe that inequality
\begin{equation}\label{puf}
\Big|x |x|^{\l-2}+|x|^\l \n \psi \Big|^2\geq  \frac{\delta^2\delta_0^2}{2}|x|^{2\l}, \quad \forall x\in \Omr,
\end{equation}
is verified for $\delta $ as in \eqref{bosgo}. Indeed,  due to Lemma \ref{crucial} we have
\begin{align}\label{gic1}
\Big|x |x|^{\l-2}+|x|^\l \n \psi \Big|^2&=|x|^{2\l-2}+|x|^{2\l}|\n \psi|^2+2 |x|^{2\l-2}x\cdot \n \psi\nonumber\\
&\geq |x|^{2\l}(|\n \psi|^2-2 \delta D_{\Omega, \psi_1})\nonumber\\
&\geq |x|^{2\l} (\delta^2\delta_0^2-2\delta D_{\Omega, \psi_1})\nonumber\\
&\geq \frac{\delta^2\delta_0^2}{2} |x|^{2\l}, \quad \forall x\in \Omega\setminus\overline{\omega_0},
 \end{align}
since $\delta\geq 4D_{\Omega, \psi_1}/
\delta_0^2$. Next, applying \eqref{puf}  and coupling the terms  independent of $\l$ in \eqref{gic2} we obtain
\begin{align*}
|\n \tau|^2&\geq (2|x|\psi-|x|^2|\n \psi|)^2 + \l^2\frac{\delta^2\delta_0^2}{2}\cat3 \phi^2 \nonumber\\
&+4\la \psi \frac{\phi}{r_0^\l}(|x|^{\l-1}+|x|^\l \n \psi\cdot x)
\nonumber\\
&+2\la \frac{\phi}{r_0^\l}|x|^2(x\cdot \nabla \psi|x|^{\l-2}+|x|^\l|\n \psi|^2), \quad \forall x\in \Omr.
\end{align*}
This and Lemma \ref{crucial} lead to (for $\l\geq 1$)
\begin{align}
|\n \tau|^2 & \geq \l^2 \frac{\delta^2\delta_0^2}{2} \cat3 \phi^2  - 2 \l \catt \phi \delta D_{\Omega, \psi_1} R_\Omega^2 \left(2\psi+1\right)\nonumber\\
& \geq \l^2 \frac{\delta^2\delta_0^2}{4} \cat3 \phi^2 \nonumber\\
& \geq \l^2 \cat3 \phi^2, \quad \forall x\in \Omr.
\end{align}
The last two inequalities above are valid since $\phi> \l \psi$, $\psi\geq \delta$ and $\delta \geq \max \{1, 2/\delta_0, 24D_{\Omega, \psi_1}R_\Omega^2/\delta_0^2 \}$.

Inequality \eqref{pici3} is trivial due to the $C^2$ regularity of $\tau$. With that,  we end the proof of Proposition \ref{tzuca}.
\end{proof}

\section{Proofs of lemmas from Section \ref{sec:1}}\label{desteptutz}

\begin{proof}[Proof of Lemma \ref{pozbd}]
 It suffices to prove that $\nabla \sigma \cdot \vec{n} \geq 0$ for all  $(x, t) \in \Gamma \times (0, T)$. First,  we have
\begin{equation}\label{gradietsigma}
\nabla \sigma=\theta(t)\Big(- 2 x\psi - |x|^2 \n \psi-\caat (\lambda
x|x|^{\lambda-2}+\lambda|x|^{\lambda} \nabla \psi)\phi \Big).
\end{equation}
The first two conditions in  \eqref{equ1}
 yield $$\nabla \psi \cdot \vec{n}=-|\nabla \psi|,\q  \forall
x\in \Gamma,$$
and
$$\nabla \sigma \cdot
\vec{n}=\theta(t)\Big(-2 x\cdot \vec{n} + |x|^2 |\n \psi|+\lambda
\caat\phi |x|^{\lambda-2}(|x|^2|\n \psi|- x\cdot \vec{n})\Big), \quad \forall x\in \Gamma.$$
 In consequence,  \eqref{impcond1} implies
$$\n \sigma \cdot \vec{n}\geq \theta(t) \left(
|x|^2(|\n \psi|-2C_\Omega)+ \l \phi \catt (|\n
\psi|-C_\Omega)\right), \quad \forall x\in \Gamma,$$
 which is nonnegative
 since $\psi$ satisfies the third condition in \eqref{equ1}. This completes the proof of Lemma \ref{pozbd}.
\end{proof}

\begin{proof}[Proof of Lemma \ref{lemaIll}]
\noindent{\textit{Computations for $I_l$.}}\\
\noindent  Next we split $I_l$ in two parts as $I_l=I_l^1+I_l^2$
where
\begin{align}
I_l^1&= -2 s \iiQ D^2 \sigma (\n z, \n z) \dx\dt - s \iiQ \D \sigma
\alpha |\n z|^2 \dx\dt + 2 \mu s \iiQ \frac{z^2}{|x|^4}x \cdot
\n \sigma \dx \dt , \nonumber\\
I_l^2&=\frac{s}{2}\iiQ z^2 \Delta^2 \sigma \(1+\alpha\) \dx \dt
        + s \iiQ z^2 \nabla \alpha \cdot \nabla \Delta\sigma \dx \dt
        +\nonumber\\
         &+\frac{s}{2} \iiQ z^2 \Delta\sigma \ \Delta \alpha \dx
         \dt+ \mu s \iiQ \frac{z^2}{|x|^2} \Delta \sigma  \
        \alpha \dx \dt.
\end{align}

Moreover, we split $I_l^1$ as $I_l^{1}=I_{l, x^2}^1+I_{l, \phi}^1$
where
\begin{align}
I_{l, x^2}^1&= -2 s \iiQ D^2 \sigma_{x^2} (\n z, \n z) \dx\dt - s
\iiQ \D \sigma_{x^2} \alpha |\n z|^2 \dx\dt\nonumber\\
& + 2 \mu s \iiQ
\frac{z^2}{|x|^4}x \cdot \n \sigma_{x^2} \dx \dt,
\end{align}
\begin{align}
I_{l, \phi}^1&=-2 s \iiQ D^2 \sigma_\phi (\n z, \n z) \dx\dt - s
\iiQ \D \sigma_{\phi} \alpha |\n z|^2 \dx\dt \nonumber\\
&+ 2 \mu s \iiQ
\frac{z^2}{|x|^4}x \cdot \n \sigma_{\phi} \dx \dt.
\end{align}
\noindent \textit{ Estimates for $I_{l, x^2}^{1}$:}

 Using the
relations \eqref{catevaformule}-\eqref{id1:2} we have
\begin{multline}\label{chica}
I_{l, x^2}^1= 4s \iiQ \theta \left( |\n z|^2 -\mu
\frac{z^2}{|x|^2}\right)\psi \dx \dt +  8 s \iiQ \theta  (x\cdot \n
z)(\n \psi \cdot \n z) \dx \dt \\
    + 2s \iiQ \theta |x|^2 D^2 \psi (\n z, \n
z)\dx \dt- 2\mu s \iiQ \theta \frac{z^2}{|x|^2} (\n \psi \cdot x)
\dx \dt - s \iiQ \D \sigma_{x^2} \alpha |\n z|^2 \dx \dt.
\end{multline}

Next, we estimate the first term in $I_{l, x^2}^1$ applying  the result of Proposition \ref{prop2}. 
Firstly, by integration by parts we
get the identities \be\label{int2} \into z \n z\cdot \n \psi \dx = -
\frac{1}{2} \into z^2 \D \psi \dx,
  \ee \be\label{int1}
  \into |x|^{2-\gamma} z \n z\cdot \n \psi \dx = - \frac{1}{2} \into
   |x|^{2-\gamma} \D \psi  z^2 \dx - \frac{(2-\gamma)}{2} \into (x\cdot \n
  \psi) |x|^{-\gamma} z^2 \dx.
 \ee
 Secondly, we apply inequality \eqref{imprhardy} for
 $u:=z\sqrt{\psi}$. Then, after integrating in time, inequality \eqref{imprhardy} becomes
\begin{multline}\label{impl1}
C_2 \iiQ \theta \psi z^2  \dx \dt + \iiQ \theta \left(|\n z|^2 -\mu
\frac{z^2}{|x|^2} \right)\psi \dx \dt\\
 + \frac{1}{4} \iiQ \theta z^2
\frac{|\n \psi|^2}{\psi} \dx \dt
 + \iiQ z \n z \cdot \n \psi \dx\dt\\
 \geq  C_3  \iiQ \theta \left( |x|^{2-\gamma} |\n z|^2 +
\frac{z^2}{|x|^\gamma}\right) \psi \dx \dt + \frac{C_3}{4} \iiQ \theta |x|^{2-\gamma}
\frac{|\n \psi|^2}{\psi} z^2 \dx \dt\\
+C_3 \iiQ \theta  |x|^{2-\gamma}  z \n z \cdot \n \psi  \dx \dt.
\end{multline}
  According to \eqref{impl1}, \eqref{int2} and \eqref{int1} 
 we obtain
\begin{multline}\label{imp1}
C_2 \iiQ \theta \psi z^2  \dx \dt + \iiQ \theta \left(|\n z|^2 -\mu
\frac{z^2}{|x|^2} \right)\psi \dx \dt\\
 + \frac{1}{4} \iiQ \theta z^2
\frac{|\n \psi|^2}{\psi} \dx \dt
 - \frac{1}{2} \iiQ \theta z^2 \D \psi \dx\dt\\
 \geq  C_3  \iiQ \theta \left( |x|^{2-\gamma} |\n z|^2 +
\frac{z^2}{|x|^\gamma}\right) \psi \dx \dt + \frac{C_3}{4} \iiQ \theta
|x|^{2-\gamma} \frac{|\n \psi|^2}{\psi} z^2 \dx \dt \\
- \frac{C_3}{2} \iiQ  \theta |x|^{2-\gamma} \D \psi z^2 \dx \dt -
C_3 \frac{(2-\gamma)}{2}\iiQ \theta (x\cdot \n \psi ) |x|^{-\gamma}
z^2 \dx \dt.
\end{multline}
Let us now compare the singular terms on the right hand side of \eqref{imp1} as follows. From the election of $r_0$ in \eqref{condr0} it is verified   $$\frac{C_3\psi}{2|x|^\gamma}\geq \frac{C_3}{2} (2-\gamma) |D\psi|_\infty
|x|^{1-\gamma}, \quad \forall x\in \tilde{\Omega}_{r_0},$$ 
and  we obtain
\begin{align}\label{prost}
\frac{C_3}{2}  \iiQ \theta \left( |x|^{2-\gamma} |\n z|^2 +
\frac{z^2}{|x|^{\gamma}} \right) & \psi \dx \dt - C_3 \frac{(2-\gamma)}{2}\iiQ \theta (x\cdot \n \psi ) |x|^{-\gamma}
z^2 \dx \dt\nonumber\\
& \geq -C_3 \frac{(2-\gamma)}{2} \sup_{r_0\leq |x|\leq R_\Omega}\{|x|^{1-\gamma}\} |D \psi|_{\infty}  \iiQs \theta z^2 \dx\dt.
\end{align}
Combining \eqref{prost} and \eqref{imp1} we deduce
\begin{multline}\label{ciucu}
\iiQ \theta \left(|\n z|^2 -\mu \frac{z^2}{|x|^2} \right) \psi \dx \dt
\geq \frac{C_3}{2} \iiQ \theta \left( |x|^{2-\gamma} |\n z|^2 +
\frac{z^2}{|x|^\gamma} \right)\psi\dx \dt \\
- \bigg( C_2 |\psi|_{\infty} + \frac{|D
\psi|_{\infty}^2}{4}+\frac{1}{2} |D^2 \psi|_{\infty} +\frac{C_3}{2}
R_\Omega^{2-\gamma} |D^2 \psi|_{\infty}+\\
+\frac{C_3}{2}(2-\gamma) \sup_{r_0\leq |x|\leq R_\Omega}\{|x|^{1-\gamma}\} |D\psi|_\infty\bigg)\iiQ
\theta z^2 \dx\dt.
\end{multline}
Reconsidering the constants in \eqref{ciucu}, there exists $ C_4$ depending on $\Omega$,  $\psi$ and $\gamma$ such that
\begin{multline}\label{condo1}
 \iiQ \theta \left( |\n z|^2 -\mu\frac{z^2}{|x|^2} \right)\psi \dx\dt
\geq \frac{C_3}{2}  \iiQ \theta \left(|x|^{2-\gamma} |\n z|^2
+\frac{z^2}{|x|^\gamma}\right) \psi \dx\dt \\
- C_4  \iiQ \theta z^2 \dx\dt.
\end{multline}
From \eqref{condo1} and \eqref{chica} we obtain
\begin{multline}\label{eqq2}
I_{l, x^2}^1 \geq 2s C_3 \iiQ \theta \left( |x|^{2-\gamma}|\n
z|^2 + \frac{z^2}{|x|^\gamma}\right)\dx \dt - 4 s C_4   \iiQ \theta
z^2 \dx \dt
\\
 -8 s |D\psi|_{\infty} \iiQ \theta |x||\n z|^2 \dx \dt - 2 s
|D^2 \psi|_{\infty} \iiQ \theta |x|^2 |\n z|^2 \dx \dt \\
-2 s |\mu| |D\psi|_{\infty} \iiQ \theta \frac{z^2}{|x|} \dx \dt -
s \iiQ \D \sigma_{x^2} \alpha |\n z|^2 \dx \dt.
\end{multline}%

Since $\gamma>1$,  the terms in \eqref{eqq2} involving the quantities $|x|^{2-\gamma}|\n z|^2, z^2/|x|^{\gamma}$ dominate the terms involving $|x||\n z|^2, |x|^2|\n z|^2$ respectively $z^2/|x|$ close enough to the origin (more precisely, for $x\in \tilde{\Omega}_{r_0}$). This is true due to
$$C_3|x|^{2-\gamma} \geq 8  |D\psi|_\infty|x|+ 2
|D^2 \psi|_\infty |x|^2 \textrm{ and }  \frac{C_3}{ |x|^\gamma} \geq \frac{|\mu|
|D\psi|_{\infty}}{|x|}, \quad  \forall x\in \tilde{\Omega}_{r_0},$$
 from the election of  $r_0$ in \eqref{condr0}. Hence,  from \eqref{eqq2} we easily
obtain
\begin{multline}\label{eqq5}
I_{l, x^2}^{1} \geq  s C_3 \iiQ \theta \left( |x|^{2-\gamma}|\n
z|^2 +
\frac{z^2}{|x|^\gamma}\right)\dx \dt- C_5 s \iiQ \theta z^2 \dx\dt \\
 - s \iiQ \D \sigma_{x^2}
\alpha |\n z|^2  \dx \dt- C_6 s \iiQs \theta |\n z|^2 \dx\dt,
\end{multline}
for some constants $C_5$, $C_6$ depending on $\Omega$,  $\psi$ and $\mu$.

\noindent \textit{Estimates for $I_{l, \sigma_\phi}^1$:}

In order to get rid of the gradient terms with negative sign in
\eqref{eqq5} we have to estimate from below the quantity
$\mathcal{T}:=I_{l, \sigma_\phi}^{1} -s \iiQ \D \sigma_{x^2} \alpha
|\n z|^2 \dx \dt-C_6 s \iiQs \theta |\n z|^2 \dx\dt$,  that is
\begin{align}\label{tavi}
\mathcal{T}&= -2 s \iiQ D^2 \sigma_\phi (\n z, \n z) \dx\dt - s
\iiQ \D \sigma_{\phi} \alpha |\n z|^2 \dx\dt + 2 \mu s \iiQ
\frac{z^2}{|x|^4}x \cdot \n \sigma_{\phi} \dx \dt \nonumber\\
&-s \iiQ \D \sigma_{x^2} \alpha |\n z|^2 \dx\dt-C_6 s \iiQs \theta |\n z|^2 \dx\dt.
\end{align}

To do that, according to Propositions \ref{prop3}-\ref{prop4} we remark that
\begin{align}
2 D^2 \tau_\phi (\n z, \n z)+\D \tau_{\phi} \alpha |\n z|^2 +\D \tau_{x^2}\alpha  |\n z|^2 &\geq \l
\left(\frac{|x|}{r_0}\right)^{\l-2}\phi |\n z|^2,\q \forall  x\in \tilde{\Omega}_{r_0},\label{pic1}\\
\big|2 D^2 \tau_\phi (\n z, \n z)+\D \tau_{\phi} \alpha |\n z|^2 +(\D \tau_{x^2}\alpha-C_6)  |\n z|^2\big|&\leq
C_7 \l^2 \catt \phi |\n z|^2, \q \forall  x\in \omega_0, \label{pic2}
\end{align}
\begin{align}
2 D^2 \tau_\phi (\n z, \n z)+\D \tau_{\phi}\alpha |\n z|^2+(\D \tau_{x^2}\alpha-C_6)  |\n z|^2 &\geq C_8\l^2
\catt \phi |\n z|^2, \q \forall  x\in \Omr,\label{pic3}
\end{align}
for $\l$ large enough and some positive constants $C_7, C_8$ uniform in $\l$. On
the other hand, it holds that
\begin{equation}\label{condo}
\frac{2 |\mu| |x\cdot \n \tau_\phi|}{|x|^4}\leq C_9 \l
\left(\frac{|x|}{r_0}\right)^{\la-4} \phi, \q \forall x\in \Omega.
\end{equation}
for some constant $C_9>0$.
 Therefore, it follows from above that
\begin{multline}\label{eqq4}
\mathcal{T}\geq C_8 s\l^2 \iiQr \theta \catt \phi |\n z|^2 \dx \dt+ \\
+ s\l \iir0 \theta \left(\frac{|x|}{r_0}\right)^{\l-2}|\n z|^2 \dx
\dt -C_9 s \l \iiQ \theta \left(\frac{|x|}{r_0}\right)^{\l-4}\phi z^2 \dx \dt\\
- C_7 s \l^2 \iiq \theta \catt \phi |\n z|^2 \dx \dt
\end{multline}
 Summing the terms in
\eqref{eqq5} and \eqref{eqq4} we get
\begin{multline}
I_l^1 \geq C_{3} s  \iiQ \theta \left( |x|^{2-\gamma}|\n z|^2 +
\frac{z^2}{|x|^\gamma}\right)\dx \dt +C_8 s\l^2 \iiQr \theta \catt \phi |\n z|^2 \dx \dt \\
+ s\la \iir0 \theta \left(\frac{|x|}{r_0}\right)^{\l-2}|\n z|^2 \dx
\dt -C_9 s \l \iiQ \theta \left(\frac{|x|}{r_0}\right)^{\l-4}\phi z^2 \dx \dt\\
- C_5 s \iiQ \theta z^2 \dx\dt - C_7 s\la^2 \iiq \theta \catt \phi
|\n z|^2 \dx \dt,
\end{multline}

 \noindent \textit{Estimates for $I_l^2$.}

  Making
use of the support of $\alpha$ located far from the origin and the $C^4$ regularity of $\tau_\phi$ we note that
\begin{align}\label{form9}
|\D^2 \tau_\phi|, \q |\D \tau_\phi|, \q |\n \D \tau_\phi|,\q
\left|\alpha \frac{\D \tau_\phi}{|x|^2}\right|  &\leq A_{\la}, \q \forall x\in \Omega,
\end{align}
where $A_{\la}$ is a big enough constant  ($A_\la$ also depends on $\Omega$,  $\psi$ and $r_0$,  but we need to emphasize its dependence on $\la$).   Then we
get
\begin{equation}
I_l^2 \geq -A_{\la} s  \iiQ \theta  |z|^2 \dx \dt.
\end{equation}
 Next we conclude
\begin{multline}
I_l \geq C_{3} s  \iiQ \theta \left( |x|^{2-\gamma}|\n z|^2 +
\frac{z^2}{|x|^\gamma}\right)\dx \dt +C_8 s\l^2 \iiQr \theta \catt \phi |\n z|^2 \dx \dt\\
+ s\la \iir0 \theta \left(\frac{|x|}{r_0}\right)^{\l-2}|\n z|^2 \dx
\dt 
- B_{\la} s \iiQ \theta z^2 \dx\dt\\
 - C_7 s\la^2 \iiq \theta \catt
\phi |\n z|^2 \dx \dt,
\end{multline}
where $B_\la=C_5+A_{\la}+C_9 \la \sup_{x\in \Omega} \{(|x|/r_0)^{\la-4}
\phi\}$.
\end{proof}

\begin{proof}[Proof of Lemma \ref{Inlfinal}]
We split  $I_{nl}=I_{nl, 1}+ I_{nl, 2}$, where $I_{nl, 1}$ are the
integrals in $I_{nl}$ restricted to $\tilde{\Omega}_{r_0}$ and $I_{nl, 2}$
are the terms in $I_{nl}$ restricted to $\Oms$.
We put $\sigma=-\theta \tau$. Then $I_{nl}$ could be written as
\begin{multline}
I_{nl}= 2 s^3\iiQ \theta^3 z^2 D^2\tau (\n \tau, \n \tau)\dx \dt
- s^3 \iiQ \theta^3 z^2 \alpha \D \tau |\n \tau|^2\dx \dt\\
-\frac{s^2}{2}\iiQ \theta^2 \alpha^2 z^2 |\D \tau|^2\dx \dt.
\end{multline}

\noindent {\textit{Computations for $I_{nl, 1}$.}}
From \eqref{bambu1} and \eqref{pici} we have that
\be\label{chicci}
T_2\geq -\l^2 \catt \phi |D \psi|_{\infty}^{2} (8 \psi^2 +2)|x|^{2}, \quad \forall x\in \tilde{\Omega}_{r_0},
\ee
which combined with \eqref{tumpi1} leads to
$T_2+T_3\geq 0$ in $\tilde{\Omega}_{r_0}$,
since $r_0\leq 1/ (|D\psi|_{\infty} \sqrt{8|\psi|_{\infty}^{2}+2})$. In consequence, due to \eqref{chici1} we
 obtain
\begin{align}\label{urs3}
2D^2\tau (\n \tau, \n \tau )- \alpha \D \tau |\n \tau|^2\geq |\n
\tau|^2, \q \forall x\in
\tilde{\Omega}_{r_0},
\end{align}
Again, from \eqref{chici1} in Proposition \ref{tzuca} we get
  \begin{align}\label{ursulet2}
2D^2\tau (\n \tau, \n \tau )- \alpha \D \tau |\n \tau|^2 &\geq |x|^2, \q \forall x\in \tilde{\Omega}_{r_0},
\end{align}
and since $\alpha$ is supported away from the origin we also have
\begin{align}\label{ursulet}
\alpha^2 |\D \tau|^2 &\leq C_{\l} |x|^2, \q \forall x\in \tilde{\Omega}_{r_0},
  \end{align}
  for some constant $C_{\l}$ depending on $\l$.
Combining \eqref{ursulet} with \eqref{ursulet2} and \eqref{urs3}, there exists $s_0=s_0(\l)$ large enough such that for any $s\geq s_0$
\begin{equation}\label{eqq7}
I_{nl,1} \geq \frac{s^3}{2}\iir0 \theta^3 |x|^2 z^2 \dx \dt.
\end{equation}

\noindent {\textit{Computations for $I_{nl, 2}$.}}

According to Propositions \ref{prop8}-\ref{prop6} and \eqref{pici2}
observe that
\begin{align}
2D^2\tau (\n \tau, \n \tau )- \alpha \D \tau |\n \tau|^2&\geq C_{10}
\l^2
\catt \phi |\n \tau|^2,\nonumber\\
& \geq C_{11}\l^4
\left(\frac{|x|}{r_0}\right)^{3\l}\phi^3, \q \forall x\in \Omr.\label{bau1}
\end{align}
In addition, it holds
\begin{align}\label{tic}
\alpha^2 |\D \tau|^2 &\leq C_{12} \l^4 \left(\frac{|x|}{r_0}\right)^{2\l}
\phi^2, \q \forall  x\in \Oms, \nonumber\\
\big|2D^2\tau (\n \tau, \n \tau )- \alpha \D \tau |\n \tau|^2\big|
&\leq C_{13} \l^2 \catt \phi |\n \tau|^2, \nonumber\\
&\leq C_{14}\l^4 \left(\frac{|x|}{r_0}\right)^{3\l}\phi^3, \q
\forall x\in \omega_0.
\end{align}
for some constants $C_{10}, C_{11}, C_{12}, C_{13}, C_{14}$ uniform in $\l$. The first two inequalities in \eqref{tic}  are consequences  of \eqref{tavico1}-\eqref{tavico2}, whereas the latter one follows from \eqref{pici3} in Proposition \ref{tzuca}.

 Then we obtain
\begin{multline}\label{eqq6}
I_{nl,2} \geq C_{11}  s^3 \l^4 \iiQr \theta^3
\left(\frac{|x|}{r_0}\right)^{3\l} \phi^3  z^2 \dx \dt- C_{14}
s^3\l^4 \iiq \theta^3
\left(\frac{|x|}{r_0}\right)^{3\l}\phi^3 z^2 \dx\dt\\
-C_{12} s^2\l^4 \iiQs \theta^2 \left(\frac{|x|}{r_0}\right)^{2\l}  \phi^2 z^2 \dx \dt.
\end{multline}
By summing the terms in \eqref{eqq6}-\eqref{eqq7}, there exists $s_0=s_0(\l)$ large enough such that for any $s> s_0$ it is satisfied
\begin{multline}
I_{nl}\geq \frac{s^3}{2} \iir0 \theta^3 |x|^2 z^2 \dx \dt+ C_{15}
s^3 \l^4 \iiQr \theta^3 \left(\frac{|x|}{r_0}\right)^{3\l} \phi^3
z^2 \dx \dt\\
- C_{16} s^3\l^4 \iiq \theta^3 \left(\frac{|x|}{r_0}\right)^{3\l}\phi^3
  z^2 \dx\dt,
\end{multline}
where $C_{15}=C_{11}/2$ and $C_{16}=C_{12}+C_{14}$.
\end{proof}

\begin{proof}[Proof of Lemma \ref{Ir}]
 According
to the expression of $\theta$ we obtain
$$|\theta'|\leq C \theta^{1+1/k}, \q |\theta''|\leq C \theta^{1+2/k},$$
for some positive constant $C$.
On the other hand, from the definition of $\sigma$ we get %
\begin{align}\label{goodbounds}
|\D \sigma|&\leq D_{\l}   \theta , \q \forall  x\in
\Omega,\nonumber\\
|\p_t \sigma|&\leq D_{\l}  \theta', \q \forall  x\in
\Omega,\nonumber\\
\p_t(|\n \sigma|^2)&\leq D_{\l} \theta \theta' |x|^2, \q \forall
x\in \tilde{\Omega}_{r_0},
\nonumber\\
\p_t(|\n \sigma|^2)&\leq D_{\l} \theta \theta'\left(\frac{|x|}{r_0}\right)^{2\la}\phi^2, \quad
\forall  x\in \Oms,
\end{align}
for a big enough  constant $D_{\la}>0$ depending on $\la$. Since
$\alpha$ is supported far from the origin we can write
\begin{multline}\label{fin1}
s^2 \iiQ \Big| \alpha z^2 \Delta \sigma \p_t \sigma \Big| \dx\dt
\leq \frac{4D_{\l}^2}{r_0^2} s^2
 \iir0 \theta^{2+1/k} |x|^2z^2 \dx \dt \\
 +
D_{\l}^2 s^2   \iiQs \theta^{2+ 1/k}  z^2\dx \dt.
\end{multline}
From \eqref{goodbounds} we also obtain
\begin{multline}\label{fin2}
 s^2  \Big| \iiQ   z^2  \p_t \left(|\n \sigma|^2\right)  \dx\dt \Big| \leq D_{\l} s^2
\iir0 \theta^{2+1/k} |x|^2z^2 \dx \dt \\
 + D_{\l} s^2
 \iiQs \theta^{2+ 1/k} \left(\frac{|x|}{r_0}\right)^{2\l} \phi z^2\dx \dt.
\end{multline}
Now,  we put
\begin{equation*}
R:= -\frac{s}{2}\iiQ z^2 \p_{tt}^{2} \sigma \dx \dt -B_\la s \iiQ
\theta z^2 \dx\dt,
\end{equation*}
where $B_\la$ is chosen as in Lemma \ref{lemaIll}.

 Then, we remark that there exists $E_\la>0$  such that
\begin{align}
\Big|B_\la s \iiQ \theta z^2 \dx\dt\Big|, \ \Big|\frac{s}{2}\iiQ z^2 \p^2_{tt} \sigma \dx \dt \Big| & \leq
E_{\l}
s\iiQ \theta^{1+ 2/k} z^2 \dx \dt.
\end{align}
Summing up these bounds we obtain
\begin{equation}\label{R}
|R|\leq  2s E_{\l}   \iiQ \theta^{1+2/k} z^2 \dx \dt.
\end{equation}
Next we write
\begin{align*}
\iiQ \theta^{1+2/k} |z|^2 \dx \dt = \iiQ
\left(\beta\theta^{1+2/k-1/q'}|x|^{\gamma/q'}|z|^{2/q} \right)
\left(\frac{1}{\beta}
\theta^{1/q'}|x|^{-\gamma/q'}|z|^{2/q'}\right)\dx\dt.
\end{align*}
Now let us take
$$q= \frac{2+\gamma}{\gamma}, \quad q'= \frac{\gamma+2}{2}.$$ Note that
$1/q+1/q'=1$, and applying the Young inequality we obtain
\begin{align}\label{R2}
\iiQ \theta^{1+2/k} z^2 \dx \dt &\leq \frac{\beta^q}{q} \iiQ
\theta^{(1+2/k-1/q')q}|x|^2 z^2 \dx \dt+ \frac{1}{q'\beta^{q'}}\iiQ
\theta
\frac{z^2}{|x|^\gamma}\dx \dt\nonumber\\
&= \frac{\beta^{q}}{q} \iiQ \theta^3|x|^2 z^2 \dx \dt +
\frac{1}{q'\beta^{q'}}\iiQ \theta \frac{z^2}{|x|^\gamma}\dx \dt,
\end{align}
provided $k=1+2/\gamma$ and $\beta>0$. Therefore, from \eqref{R} and \eqref{R2} we
have
\begin{equation}\label{fin3}
|R|\leq 2 E_\la s \left(\beta^q \iiQ \theta^3 |x|^2 z^2\dx \dt +
\frac{1}{\beta^{q'}}\iiQ \theta \frac{z^2}{|x|^\gamma}\dx
\dt\right).
\end{equation}
 Consequently, from \eqref{fin1}-\eqref{fin3} it follows that
\begin{multline}
|I_r|\leq F_{\l}\Bigg(s^2 \iir0 \theta^3 |x|^2 z^2 \dx \dt+ s
\beta^q \iiQ \theta^3 |x|^2z^2\dx \dt\\
+\frac{s}{\beta^{q'}} \iiQ \theta \frac{z^2}{|x|^\gamma} \dx\dt +
s^2 \iiQs \theta^3 \left(\frac{|x|}{r_0}\right)^{3\l}z^2\dx \dt\Bigg),
\end{multline}
for a new constant $F_\la>0$.

 Take $\beta$ such that $F_\la/\beta^{q'}=C_{3}/2$. Then
there exists $s_0(\l)$ such that for $s\geq s_0(\l)$ we finish
the proof of Lemma \ref{Ir}.
\end{proof}
\noindent {\bf Acknowledgements.}
The author would like to thank    Sylvain Ervedoza, Enrique Zuazua
and Jean-Pierre Puel for useful suggestions related to this paper.

This work was supported by the both grants of the Ministry of National Education, CNCS  UEFISCDI Romania, project PN-II-ID-PCE-2012-4-0021  and project  PN-II-ID-PCE-2011-3-0075, and the Grant  MTM2011-29306-C02-00 of the MICINN (Spain).

This work started when the author was PhD Student within the doctoral program of Universidad Aut\'{o}noma de Madrid and member in the PDE research line at  BCAM - Basque Center for Applied Mathematics, in Bilbao. The author would like to thank both institutions for their support.

\section{Appendix: Proof of Proposition \ref{prop2}}\label{7sec}

\begin{proof}[Proof of Proposition \ref{prop2}]
Let us check the validity of Proposition \ref{prop2}. In view of Proposition \ref{prop1}, it suffices to prove  the inequality
  \be\label{cruineq}
\forall w\in \hoi, \quad \into|x|^{2-\gamma}|\nabla w|^2\dx\leq
R_{\Omega}^{2-\gamma}\left(\into|\nabla
w|^2\dx-\frac{N^2}{4}\into\frac{w^2}{|x|^2}\dx\right)
+C\into
w^2\dx,
  \ee
  for some constant $C$ depending on $\Omega$ and $\gamma$.

  In the sequel we generalize the proof given in \cite{cristisumJFA} for $\gamma=0$ and extend it to any $\gamma\in [0, 2)$. Here we reproduce the main steps of the proof in \cite{cristisumJFA} adapted to the case $\gamma\in [0, 2)$.\\

\noindent {\bf Step 1.} Firstly we show that inequality \eqref{cruineq} is true  in a neighborhood of $x=0$. More precisely, there exists   $r_1>0$ small enough depending on $\Omega$ and there exists $C\in \rr$  depending on $\Omega$ and $r_1$ such that
\begin{equation}\label{eqr0}
\int_{\Omg}|x|^{2-\gamma}|\nabla w|^2\dx\leq
R_{\Omega}^{2-\gamma}\Big(\int_{\Omg}|\nabla
w|^2\dx-\frac{N^2}{4}\int_{\Omg}\frac{w^2}{|x|^2}\dx\Big)
+C\int_{\Omg}
w^2\dx,
\end{equation}
holds true for any function $w\in C_{0}^{\infty}(\Omg)$,
where $\Omg=\Omega \cap B(0, r_1)$.

Since the result in \eqref{eqr0} is local, without losing generality it is enough to check the validity of \eqref{eqr0} for domains of the form
\begin{equation}\label{model}
\beta <0, \qquad  \Upsilon(\beta, r_1):= \{x\in \rr^N \ | \  x_N \geq \beta |x'|^2, |x|\leq r_1  \}.
 \end{equation}
Indeed, this is true since for any $\Omg$, with $r_1$ small enough, there exists $\beta<0$ depending on $r_1$ such that
$\Omg\subset \Upsilon(\beta, r_1)$.

For those reasons, the result \eqref{eqr0} valid for $ \Upsilon(\beta, r_1)$  still holds true for $\Omg$    since we can prove it
for test functions  extended from zero up to the domain $ \Upsilon(\beta, r_1)$.

Next we check the validity of Step 1 for such domains $\Omg$ as in \eqref{model}. In view of that, let us consider a
smooth function $\phi$ which satisfies
\be\label{supersol}
-\D \phi\geq \frac{N^2}{4}\frac{\phi}{|x|^2}, \q \phi >0,  \q \forall x\in \Omg.\ee
In view of \cite{Fall},  it is not difficult to see that
 \be\label{superss}
 \phi(x)=
\rho(x)e^{(1-N)\rho(x)}\left|\log \frac{1}{|x|}\right|^{1/2}|x|^{-N/2},
\ee
satisfies \eqref{supersol} for $r_1$ small enough. \\

With the transformation $w=\phi u$  for such $\phi $ as in \eqref{superss} we get
\be\label{iden}
|\n w|^2 = |\n \phi|^2 u^2 +\phi^2 |\n u|^2 +2 \phi u \n \phi \cdot \n u.
\ee
Integrating we obtain
\begin{align}\label{lema2}
\intos |\n w|^2\dx &
= \intos |\n u|^2 \phi^2 \dx-\intos \frac{\D \phi}{\phi} w^2\dx.
\end{align}
On the other hand, multiplying in \eqref{iden} by $|x|^{2-\gamma}$ and integrating  we have
\begin{align}\label{eqp1}
\intos |x|^{2-\gamma} |\n w|^2 \dx&=\intos |x|^{2-\gamma} |\n \phi|^2 u^2 \dx +\intos
|x|^{2-\gamma} \phi^2 |\n u|^2 \dx\nonumber\\
&+\frac{1}{2}\intos |x|^{2-\gamma} \n (\phi^2) \cdot
\n (u^2) \dx
\end{align}
For the last term in \eqref{eqp1}   we deduce
\begin{align}\label{eqp2}
\frac{1}{2}\intos |x|^{2-\gamma} \n(\phi^2)\cdot \n (u^2) \dx&
=-(2-\gamma)\intos  |x|^{-\gamma}\frac{x\cdot \n \phi}{\phi} w^2\dx -\intos |x|^{2-\gamma} |\n
\phi|^2u^2 \dx\nonumber\\
&-\intos \frac{ \D \phi}{\phi} |x|^{2-\gamma} w^2 \dx .
\end{align}
According to \eqref{eqp1} and \eqref{eqp2} we obtain
\begin{align}\label{eqp3}
\intos |x|^{2-\gamma} |\n w|^2 \dx&= \intos |x|^{2-\gamma} \phi^2 |\n u|^2 \dx -(2-\gamma)\intos |x|^{-\gamma}
\frac{  x \cdot \n \phi}{\phi} w^2\dx \nonumber\\
&-\intos  \frac{\D \phi}{\phi}
 |x|^{2-\gamma} w^2 \dx.
\end{align}

Taking into account the election of $\phi$ in \eqref{superss} we can show that there exists $C(r_1)>0$ such that
\begin{align}\label{suma}
-\frac{\D \phi}{\phi}= \frac{N^2}{4|x|^2}+P,
\end{align}
where
\be\label{gic3}
P(x)\geq \frac{C(r_1)}{|x|^2\left(\log \frac{1}{|x|}\right)^2}, \quad \forall x\in \tilde{\Omega}_{r_1}.
\ee
 Then from \eqref{lema2} and \eqref{suma}
we have
\begin{align}\label{integr}
\intos |x|^{2-\gamma} \phi^2 |\n u|^2 \dx &\leq R_{\Omega}^{2-\gamma}\intos \phi^2 |\n
u|^2\dx= R_{\Omega}^{2-\gamma}  \Big(\intos |\n w|^2 \dx+\intos \frac{\D
\phi}{\phi}
w^2 \dx\Big)\nonumber\\
&= R_{\Omega}^{2-\gamma} \intos \Big( |\n w|^2-\frac{N^2}{4}\frac{w^2}{|x|^2}
\Big)\dx -R_{\Omega}^{2-\gamma}\intos Pw^2\dx.
\end{align}
From above and \eqref{eqp3} it follows that
\begin{align}\label{eqp4}
\intos |x|^{2-\gamma} |\n w|^2\dx &\leq R_{\Omega}^{2-\gamma} \intos \Big( |\n w|^2
-\frac{N^2}{4}\frac{w^2}{|x|^2}\Big)\dx - R_{\Omega}^{2-\gamma} \intos P
w^2\dx\nonumber\\
& -(2-\gamma) \intos |x|^{-\gamma}\frac{x\cdot \n \phi}{\phi} w^2 \dx + \intos \Big(
\frac{N^2}{4|x|^2}+P\Big)|x|^{2-\gamma} w^2 \dx\nonumber\\
& = R_{\Omega}^{2-\gamma} \intos \Big( |\n w|^2
-\frac{N^2}{4}\frac{w^2}{|x|^2}\Big)\dx +\intos(|x|^{2-\gamma}-R_{\Omega}^{2-\gamma}) P
w^2\dx\nonumber\\
&-(2-\gamma)\intos |x|^{-\gamma}\frac{x\cdot \n \phi}{\phi} w^2 \dx + \frac{N^2}{4} \intos
w^2 \dx.
\end{align}
Since $r_1\leq R_\Omega/2$, from \eqref{eqp4} and \eqref{gic3} we obtain
\begin{align}\label{gic5}
\intos |x|^{2-\gamma} |\n w|^2\dx & \leq R_{\Omega}^{2-\gamma} \intos \left( |\n w|^2
-\frac{N^2}{4}\frac{w^2}{|x|^2}\right)\dx\nonumber\\
&-\left(1- \frac{1}{2^{2-\gamma}}\right)R_{\Omega}^{2-\gamma}C(r_1) \intos \frac{w^2}{|x|^2 \left(\log \frac{1}{|x|}\right)^2} \dx\nonumber\\
& -(2-\gamma) \intos |x|^{-\gamma} \frac{x\cdot \nabla \phi}{\phi} w^2 \dx +\frac{N^2}{4} \intos \frac{w^2}{|x|^\gamma} \dx.
\end{align}
Moreover,  for $r_1\leq R_\Omega/2$ small enough we have
   $$ \q \n \rho(x)\cdot x\geq 0, \q \forall x\in \Omg.$$
 Then, we remark that
$$\frac{x\cdot \n \phi}{\phi}= \frac{x\cdot \n \rho(x)}{\rho(x)} +O(1), \qquad \textrm{ as } x\rightarrow 0,$$
 and therefore from \eqref{gic5} it follows
 \begin{align}
 \intos |x|^{2-\gamma} |\n w|^2\dx & \leq R_{\Omega}^{2-\gamma} \intos \left( |\n w|^2
-\frac{N^2}{4}\frac{w^2}{|x|^2}\right)\dx \nonumber\\
&-\left(1- \frac{1}{2^{2-\gamma}}\right)R_{\Omega}^{2-\gamma}C(r_1) \intos \frac{w^2}{|x|^2 \left(\log \frac{1}{|x|}\right)^2}\dx +C_2 \intos \frac{w^2}{|x|^\gamma} \dx,
 \end{align}
 for some constant $C_2>0$ depending on $r_1$, $R_\Omega$ and $\gamma$. Since the logarithmic term is more singular than $|x|^{-\gamma}$ as $x$ tends to zero, we obtain
 \begin{align}
 \intos |x|^{2-\gamma} |\n w|^2\dx & \leq R_{\Omega}^{2-\gamma} \intos \left( |\n w|^2
-\frac{N^2}{4}\frac{w^2}{|x|^2}\right)\dx-\frac{C_2}{r_0^\gamma} \intos w^2 \dx,
 \end{align}
 for some $r_0< r_1$.  With this we  finish the proof of Step 1.\\

\noindent {\bf Step 2.}  This step consists in applying a cut-off
argument to transfer the validity of inequality \eqref{eqr0} from
$\Omg$ to $\Omega$. More precisely, we consider a cut-off
function $\theta\in C_{0}^{\infty}(\rr^N)$ such that
\begin{equation}\theta(x)=\left\{\begin{array}{ll}
  1, & |x|\leq r_1/2, \\
  0, & |x|\geq r_1. \\
\end{array}\right.
\end{equation}
Then we split $w\in C_{0}^{\infty}(\Omega)$ as follows
\be\label{cutoff}
w=\theta w+(1-\theta) w:= w_1+w_2.
\ee

As shown in Lemma A.1 in \cite{cristisumJFA}, for a smooth function $p:
C^{\infty}(\overline{\Omega})\rightarrow \rr$ which is bounded and
non-negative,  there exists a constant $C>0$ depending on $\Omega, p, r_1$ such that the
following inequality holds
\begin{equation}\label{eq102}
\int_{\Omega}\rho(x)\nabla w_1\cdot\nabla w_2 \dx \geq
-C \int_{\Omega}|w|^2\dx.
\end{equation}
Now we are able to finalize Step 2. Indeed, splitting
$w$ as before we get
\begin{align}\label{labe}
\int_{\Omega}|x|^{2-\gamma}|\nabla w|^2\dx
&= \int_{\Omg}|x|^{2-\gamma}|\nabla
w_1|^2\dx+\int_{\Omega\setminus \tilde{\Omega}_{r_1/2}} |x|^{2-\gamma}|\nabla
w_2|^2\dx\nonumber\\
&+\int_{\Omg\setminus\tilde{\Omega}_{r_1/2}} |x|^{2-\gamma} \nabla
w_1\cdot\nabla w_2\dx
\end{align}
 Applying \eqref{eqr0} to $w_1$, from \eqref{labe}  we obtain (reconsidering the constant $C$)
\begin{align}\label{equ42}
\int_{\Omega}|x|^{2-\gamma}|\nabla w|^2\dx
 &\leq R_{\Omega}^{2-\gamma}\Big(\int_{\Omega} |\nabla w|^2\dx-\frac{N^2}{4}
 \intos \frac{w_1^2}{|x|^2}\dx\Big)+ C \int_{\Omega} w^2\dx\nonumber\\
 &-\int_{\Omg\setminus\tilde{\Omega}_{r_1/2}} 2(R_{\Omega}^{2-\gamma}-|x|^{2-\gamma})\nabla w_1\cdot\nabla
w_2\dx.
\end{align}
Considering  $\rho=2(R_{\Omega}^{2-\gamma}-|x|^{2-\gamma})$ in \eqref{eq102}, from
(\ref{equ42}) we get
\begin{equation}\label{eq104}
\into|x|^{2-\gamma}|\nabla w|^2\dx\leq R_{\Omega}^{2-\gamma}\Big(\into |\nabla
w|^2\dx-\frac{N^2}{4}\intos
\frac{w_1^2}{|x|^2}\dx\Big)+C_1\into w^2\dx,
\end{equation}
for some new constant $C_1\in \rr$.
On  the other hand we have
\begin{align}\label{equ43}
\intos\frac{w_1^2}{|x|^2} \dx
\geq \into \frac{w^2}{|x|^2}\dx-C_2\into w^2\dx,
\end{align}
for some $C_2\in \rr$.  From (\ref{eq104}) and (\ref{equ43}) the  conclusion of Proposition \ref{prop2}  yields choosing  $r_1$ small enough, $r_1< R_\Omega/2$.
\end{proof}


\bibliographystyle{plain}

\addcontentsline{toc}{section}{References}

\end{document}